\documentclass[10pt]{article}
\usepackage[dvips]{graphics}
\usepackage{epsfig,graphicx}
\usepackage{amsfonts}
\usepackage{amsmath}

\input xy
\xyoption{all}

\begin{document}

\newcommand{\lar}{\longrightarrow}
\newcommand{\kF}{\mathcal F}
\newcommand{\La}{\Lambda}
\newcommand{\la}{\lambda}
\newcommand{\de}{\delta}
\newcommand{\kQ}{\mathcal Q}
\newcommand{\kG}{\mathcal G}
\newcommand{\kL}{\mathcal L}
\newcommand{\kB}{\mathcal B}
\newcommand{\kM}{\mathcal M}
\newcommand{\kO}{\mathcal O}
\newcommand{\kP}{\mathcal P}
\newcommand{\kI}{\mathcal I}
\newcommand{\kC}{\mathcal C}
\newcommand{\kJ}{\mathcal J}
\newcommand{\bJ}{\mathbf J}
\newcommand{\bE}{\mathbf E}
\newcommand{\bF}{\mathbf F}
\newcommand{\bU}{\mathbf U}
\newcommand{\bB}{\mathbf B}
\newcommand{\kE}{\mathcal E}
\newcommand{\kT}{\mathcal T}
\newcommand{\kN}{\mathcal N}
\newcommand{\kS}{\mathcal S}
\newcommand{\kU}{\mathcal U}
\newcommand{\kk}{\mathbf k}
\newcommand{\kA}{\mathcal A}
\newcommand{\PP}{\mathbf P}
\newcommand{\bA}{\mathbf A}
\newcommand{\bI}{\mathbf I}
\newcommand{\ZZ}{\mathbb Z}
\newcommand{\RR}{R}
\newcommand{\balpha}{\mathbf \alpha}
\newcommand{\bbetha}{\mathbf \beta}
\newcommand{\FF}{\mathbf F}
\newcommand{\tens}{\otimes}
\newcommand{\Hom}{\mathop{\rm Hom}}

\renewcommand{\mod}{\mathop{\rm mod}}
\newcommand{\TF}{\mathop{\rm TF}}
\newcommand{\T}{\mathop{\rm T}}
\newcommand{\End}{\mathop{\rm End}}
\newcommand{\Ext}{\mathop{\rm Ext}\nolimits^1}
\newcommand{\TC}{\mathop{\rm TC}}

\newtheorem{theorem}{Theorem}[section]
\newtheorem{corollary}[theorem]{Corollary}
\newtheorem{remark}[theorem]{Remark}
\newtheorem{example}[theorem]{Example}
\newtheorem{lemma}[theorem]{Lemma}
\newtheorem{proposition}[theorem]{Proposition}
\newtheorem{statement}[theorem]{Statement}
\newtheorem{definition}[theorem]{Definition}
\newtheorem{conjecture}[theorem]{Conjecture}
\newtheorem{problem}[theorem]{Problem}
\newtheorem{question}[theorem]{Question}

\def\Iso{\mathop\mathrm{Iso}\nolimits}
\def\lfd{locally finite dimensional }
\def\lof{locally finite }
\def\oc{one-to-one correspondence}
\def\iff{if and only if }
\def\set#1{\left\{\,#1\,\right\}}
\def\setsuch#1#2{\left\{\,#1\,|\,#2\,\right\}}
\def\dt{_\bullet}
\def\mdl#1{{#1}\mbox{-}\rm mod}
\def\pro#1{{#1}\mbox{-}\rm pro}
\def\Mdl#1{\mathrm{Mod}\mbox{-}#1}
\def\dC{\mathfrak C}
\def\bS{\mathbf S}
\def\om{\omega}
\def\gM{\mathfrak m} 
\def\bup{\bigcup}
\def\={\setminus}         
\def\rad{\mathop\mathrm{rad}\nolimits} 
\def\End{\mathop\mathrm{End}\nolimits} 
\newcommand{\iso}{\cong}
\newcommand{\+}{\oplus}
\def\dX{\mathfrak X}
\def\dC{\mathfrak C}
\def\dG{\mathfrak G}
\def\dF{\mathfrak F}
\def\dE{\mathfrak E}
\def\ti{\tilde}
\def\8{\infty}
\def\Mod{\mathrm{Mod}\,}
\def\ob{\mathrm{Ob}\,}

\begin{center}
{\large
Derived Categories of Nodal Algebras
 \footnote{
The both authors were partially supported CRDG Grant UM2-2094. The first 
author was also partially supported by DFG Schwerpunkt ``Globale Methoden 
in der komplexen Geometrie''.}
}
\end{center}

\date{}

\begin{center}
Igor Burban and Yuriy Drozd 
\end{center}
\begin{center}
{\footnotesize Universit\"at Kaiserslautern, Fachbereich Mathematik, \\
Postfach 3049, 67653,  Kaiserslautern, Germany \\
and \\
Department of Mechanics and Mathematics,\\
Kyiv Taras Shevchenko University, 01033 Kyiv, Ukraine \\
E-mails: burban@mathematik.uni-kl.de, yuriy@drozd.org}
\end{center}


\begin{abstract}
In this article we classify  indecomposable objects of the derived categories 
of finitely-generated modules over certain infinite-dimensional algebras.
The  considered  class of algebras (which we call nodal algebras)
contains such well-known algebras as the 
complete ring of a double nodal point  $\kk[[x,y]]/(xy)$ and the completed
path algebra of 
the Gelfand quiver.  As a corollary we obtain a description 
of  the derived category
of Harish-Chandra modules over $SL_{2}({\mathbb R})$. We also give an algorithm, which  allows to construct  projective resolutions
of indecomposable complexes. In the appendix we prove the Krull-Schmidt theorem for homotopy categories. 

\emph{Key words:} derived categories; tame matrix problems; Harish-Chandra
modules; Krull-Schmidt theorem. 
\end{abstract}


\section{\bf Introduction} 
Let $A$ be a \emph{pure noetherian complete algebra}, i.e. an associative $\kk$-algebra such that:
\begin{enumerate}
\item Its center $C$ is a complete local noetherian  $\kk$-algebra. 
\item $A$ is finitely generated $C$-module without minimal submodules.
\end{enumerate}
Denote by $r$ the radical of $A$. It was shown in \cite{Pure} that
$A$ is tame if and only if it satisfies
the following conditions:

\begin{enumerate}
\item The algebra $\tilde{A} = \End_{A}(r)$ is hereditary.
\item ${\rm rad}(\tilde{A}) = r$.
\item For any simple left $A$--module $U$ the length of the left $A$-module 
$\tilde{A}\tens_{A} U$ is at most  2 and for any simple right
 $A$--module $V$ the length of the right $A$--module $V\tens_{A}\tilde{A}$
is at most  2.
\end{enumerate}

\noindent
We call  algebras satisfying these conditions \emph{nodal} algebras.  

Our description of the derived category of a nodal algebra shows that it is 
tame at least in ``pragmatic'' sense, i.e. one can obtain a list of its indecomposable  objects as a union of
 one-dimensional families and some discrete set of objects staying apart. Unfortunally, the definition of derived tameness proposed in 
\cite{GKr} can be only applied to finite dimensional algebras of finite 
global dimension and nodal algebras usually 
satisfy neither of these conditions.

The methods developed in this article can be also applied to 
finite dimensional 
gentle and skew-gentle algebras considered in \cite{Skowr}, \cite{Ringel},
\cite{laPena}, \cite{Bekkert}, \cite{Bekkert1}, as well as to some other 
algebras \cite{Burban} and to some derived categories of coherent sheaves
\cite{curves}.  An advantage of these methods is that  they  also work 
in cases, when an algebra has  infinite homological dimension and describe
the derived category of bounded from the \emph{right} complexes. 
 The developed technique allows to write down  projective resolutions of indecomposable complexes.

For the sake of simplicity we suppose that the  field $\kk$ 
is algebraically closed.
Let us rewrite  the definition of nodal algebras   in a more transparent 
 form. Let $U$ be a simple
$A$-module, $P\lar U$ its projective covering. Then we have an exact sequence
$$
  0\lar rP \lar P \lar U \lar 0.
$$
Apply the functor $\tilde{A}\tens_{A}$ to this sequence. We get
$$
  \tilde{A}\tens_{A} rP \lar \tilde{A}\tens_{A}P 
\lar \tilde{A}\tens_{A}U \lar 0.
$$
But $r = {\rm rad}(\tilde{A})$, hence 
${\rm im}(\tilde{A}\tens_{A} rP \lar \tilde{A}\tens_{A} P) = r\tens_{A} P =
{\rm rad}(\tilde{A}\tens_{A} P)$.
 So we have an exact sequence
$$
0\lar {\rm rad}(\tilde{A}\tens_{A}P) \lar \tilde{A}\tens_{A}P \lar \tilde{A}\tens_{A}U \lar 0.
$$ 
Therefore  $\tilde{A}\tens_{A}U$ is a direct sum of simple $\tilde{A}$-modules. 
Let $U_{1}, U_{2},\dots, U_{m}$ be the set of all non-isomorphic 
simple $A$-modules, 
$V_{1}, V_{2},\dots, V_{n}$ the set of all non-isomorphic simple $\tilde{A}$-modules. Consider
the graph $\Gamma$ with vertices $U_{i}, V_{j}$, $i=1,\dots,m, \, j=1,\dots, n$.  There
is an arrow from $U_{i}$ to $V_{j}$ if and only if  $V_{j}$ is a direct summand of 
$ \tilde{A}\tens_{A} U_{i}$. Then, as it was shown  in \cite{Pure},  the
last condition in the criteria of tameness is equivalent to the following
condition: all  connected components of $\Gamma$ are of the form:
\begin{enumerate}
\item
$
V' \longleftarrow U \lar V''.
$
\item $ U' \lar V \longleftarrow U''.$
\item $U\lar V. $
\end{enumerate}

\noindent
Let us consider some examples. 
\begin{example}
{\rm 
Let $A = \kk[[x,y]]/(xy)$, ${\mathfrak m} = (x,y)$ be its maximal ideal. Then
$\tilde{A} = \End_{A}({\mathfrak m}) = \kk[[x]]\times \kk[[y]]$. Let $U$ be the unique 
simple $A$-module, $V_{1}, V_{2}$ be simple $\tilde{A}$-modules. Then the 
graph 
$\Gamma$ has the form
$$
 V_{1} \longleftarrow U \lar V_{2}.
$$
}
\end{example}

\begin{example}
{\rm 
Let 
$$
A = \{
\left(
\begin{array}{cc}
f_{11} & tf_{12}  \\
f_{21} & f_{22}  \\
\end{array}
\right)| f_{ij} \in \kk[[t]], 1\le i,j \le 2; f_{11}(0) = f_{22}(0) \} 
\subseteq {\rm Mat}_{2}(\kk[[t]]). 
$$
As one can easily observe, $A$ is just the  algebra 
$\kk\langle\langle x,y\rangle\rangle/ (x^{2}, y^{2})$ 
($\kk\langle\langle x,y\rangle\rangle$ is the algebra of formal power series
in two non-commutative variables). The endomorphism algebra
of its radical is just 
$$\{
\left(
\begin{array}{cc}
f_{11} & tf_{12}  \\
f_{21} & f_{22}  \\
\end{array}
\right)| f_{ij} \in \kk[[t]], 1\le i,j \le 2 \} 
\subseteq {\rm Mat}_{2}(\kk[[t]]). 
$$
It is easy to see that it is just the completed path algebra of the quiver

\begin{figure}[ht]
\hspace{4.6cm}
\includegraphics[height=1.5cm,width=2.5cm,angle=0]{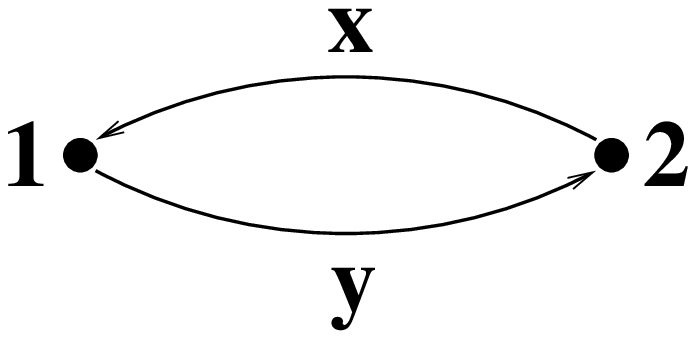}
\end{figure}

\noindent
Indeed, an  isomorphism is given by 
$$
x \mapsto 
\left(
\begin{array}{cc}
0 & t  \\
0 & 0 \\
\end{array}
\right),
y \mapsto 
\left(
\begin{array}{cc}
0 & 0 \\
1 & 0 \\
\end{array}
\right).
$$
}
\end{example}

\noindent
The graph $\Gamma$ again has the form
$$
V_{1} \longleftarrow U \lar V_{2}. 
$$
Here and further on we consider the natural completion of path algebras, namely, the $J$-adic one, where $J$ is the ideal generated by all arrows.

\begin{example}
{\rm 
Let $A$ be the completed path algebra of the Gelfand quiver

\begin{figure}[ht]
\hspace{4.2cm}
\includegraphics[height=1.5cm,width=5cm,angle=0]{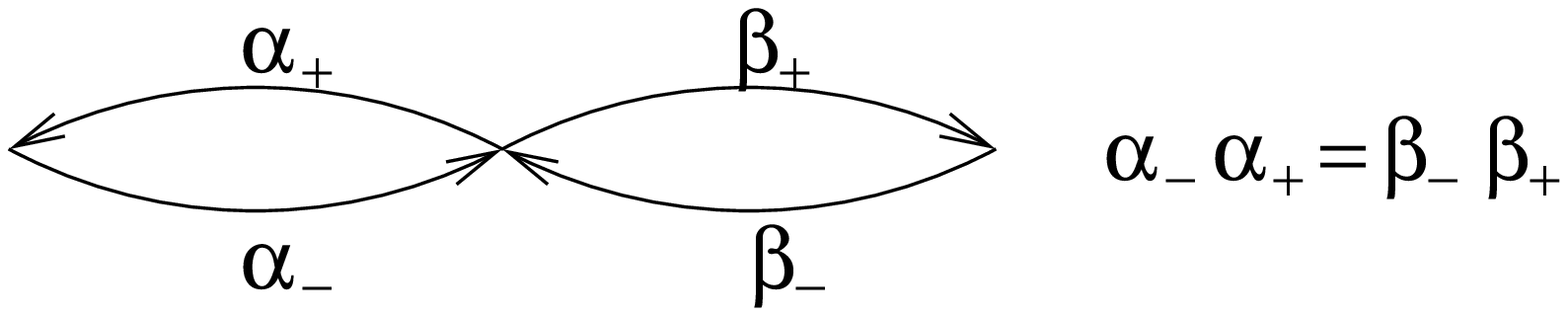}
\end{figure}

\noindent
As one can easily see, 
$$
A \cong \{
\left(
\begin{array}{ccc}
f_{11} & tf_{12} &tf_{13} \\
f_{21} & f_{22} & tf_{23} \\
f_{31} & tf_{32} & f_{33} \\
\end{array}
\right)| f_{ij} \in \kk[[t]], 1\le i,j \le 3 \} \subseteq {\rm Mat}_{3}(\kk[[t]])
$$
}
\end{example}
An  isomorphism is given by 
$$
\alpha_{+} \mapsto 
\left(
\begin{array}{ccc}
0 & 0 &0 \\
1 & 0 & 0\\
0 & 0 & 0 \\
\end{array}
\right),
\beta_{+} \mapsto 
\left(
\begin{array}{ccc}
0 & 0 &0 \\
0 & 0 & 0\\
1 & 0 & 0 \\
\end{array}
\right),
$$
$$
\alpha_{-} \mapsto 
\left(
\begin{array}{ccc}
0 & t &0 \\
0 & 0 & 0\\
0 & 0 & 0 \\
\end{array}
\right),
\beta_{-} \mapsto 
\left(
\begin{array}{ccc}
0 & 0 &t \\
0 & 0 & 0\\
0 & 0 & 0 \\
\end{array}
\right).
$$
The endomorphism algebra $\tilde{A} = \End({\rm rad}(A))$ is 
$$
\tilde{A} \cong \{
\left(
\begin{array}{ccc}
f_{11} & tf_{12} &tf_{13} \\
f_{21} & f_{22} & f_{23} \\
f_{31} & f_{32} & f_{33} \\
\end{array}
\right)| f_{ij} \in \kk[[t]], 1\le i,j \le 3 \} \subseteq {\rm Mat}_{3}(\kk[[t]]).
$$
and is Morita equivalent to 
 $$ \{ 
\left(
\begin{array}{cc}
f_{11} & tf_{12}  \\
f_{21} & f_{22}  \\
\end{array}
\right)| f_{ij} \in \kk[[t]], 1\le i,j \le 2 \} 
\subseteq {\rm Mat}_{2}(\kk[[t]]), 
$$
which is the  completed path algebra of the quiver
\begin{figure}[ht]
\hspace{4.9cm}
\includegraphics[height=1.5cm,width=2cm,angle=0]{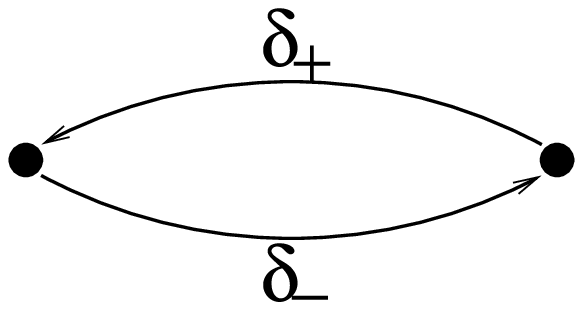}
\end{figure}

\noindent
Note that $\tilde{A}$ is isomorphic to the completed path algebra
 of the
following (non-basic) quiver

\begin{figure}[ht]
\hspace{3.5cm}
\includegraphics[height=2.2cm,width=7cm,angle=0]{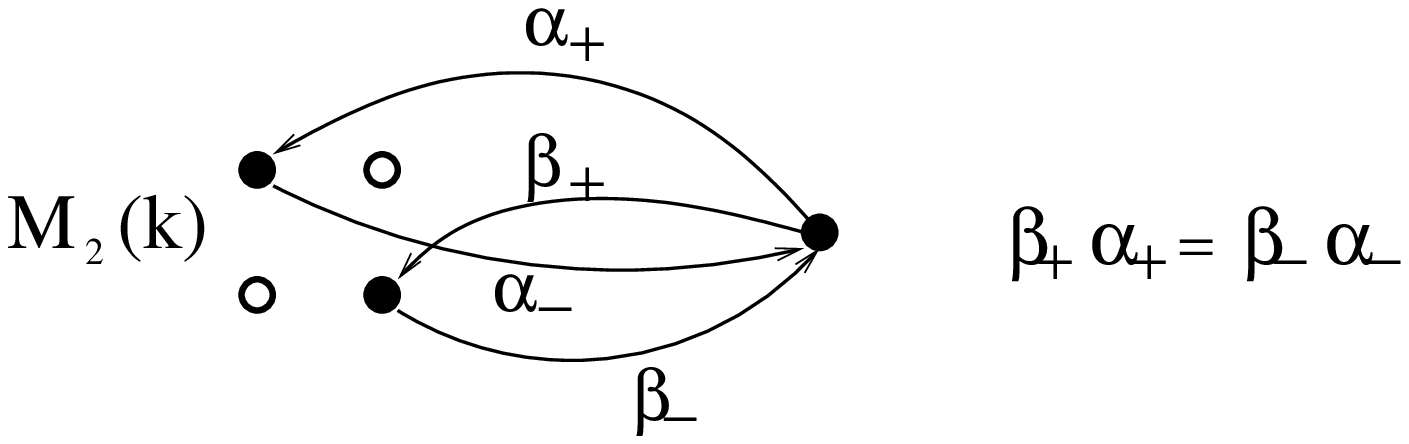}
\end{figure}

Let $U_{1}, U_{2}, U_{3}$ be simple $A$-modules, $V, W$ simple 
$\tilde{A}$-modules. Our graph $\Gamma$ has the form
$$
U_{1} \lar V \longleftarrow U_{2},\hspace{0.3cm} U_{3}\lar W. 
$$

\section{The main construction}

Let $A$ be a semi-perfect associative $\kk$-algebra 
(not necessarily finite dimensional),
$A\subset \tilde{A}$ be an embedding such that $r = {\rm rad}(A) = {\rm rad}(\tilde{A})$.
Let $I\subset A$ be a  two-sided $\tilde{A}$-ideal  containing $r$. It means 
that 
 $r\subseteq I = I\tilde{A} = \tilde{A} I$, thus
 $A/I$ and $\tilde{A}/I$ are semi-simple algebras.

Let $\tilde{A}\tens_{A}$ be the derived functor of the tensor product. 
We want to describe the fibers of the  map
$$
{\rm Ob}\bigl(D^{-}(A-\mod)\bigr) \lar 
{\rm Ob}\bigl(D^{-}(\tilde{A}-\mod)\bigr).
$$
\begin{remark}
$A-$\mbox{\rm mod} 
denotes the category of \emph{finitely-generated} $A$-modules.
We always consider objects of derived categories as complexes of projective modules.
\end{remark}

\begin{definition}
Consider the following category of triples of 
complexes $\TC_{A}$ 

\begin{enumerate}
\item Objects are  triples 
$(\widetilde {\mathcal P}_{\bullet}, {\mathcal M}_{\bullet}, i)$, where

$\tilde{\mathcal P}_{\bullet} \in D^{-}(\tilde{A}-\mod)$, 

${\mathcal M}_{\bullet} \in D^{-}(A/I-\mod)$, 

$i : {\mathcal M}_{\bullet} \lar
\tilde{A}/I\tens_{\tilde{A}} \tilde{\mathcal P}_{\bullet}$ 
 a morphism in
 $ D^{-}(A/I-\mod) $, such that 

$\tilde{i} :
\tilde{A}/I\tens_{A} {\mathcal M}_{\bullet} \lar
\tilde{A}/I\tens_{\tilde{A}}\tilde{\kP}_{\bullet}$
is an isomorphism in $D^{-}(\tilde{A}/I-\mod).$ 

\item Morphisms
 $                              
(\tilde{\mathcal P_{\bullet}}_{1}, {\mathcal M_{\bullet}}_{1}, i_{1}) \lar
(\tilde{\mathcal P_{\bullet}}_{2}, {\mathcal M_{\bullet}}_{2}, i_{2})
 $
are  pairs $(\Phi, \varphi)$, 
$$\tilde{\mathcal P_{\bullet}}_{1} \stackrel{\Phi}\lar 
\tilde{\mathcal P_{\bullet}}_{2},\quad
 {\mathcal M_{\bullet}}_{1} \stackrel{\varphi}\lar
 {\mathcal M_{\bullet}}_{2},
$$ such that 

\begin{tabular}{p{2.7cm}c}
 &
\xymatrix
{ 
\tilde{A}/I\tens_{\tilde{A}} \tilde{\mathcal P_{\bullet}}_{1} 
\ar[d]_{
\Phi\tens \mbox{\rm id}} &
{\mathcal M_{\bullet}}_{1} \ar[l]_{\qquad i_1} \ar[d]^{\varphi} \\               \tilde{A}/I\tens_{\tilde{A}} \tilde{\mathcal P_{\bullet}}_{2}  &
{\mathcal M_{\bullet}}_{2} \ar[l]^{\qquad i_2} 
}
\end{tabular}

is commutative. 
\end{enumerate}
\end{definition}

\begin{remark}
If an algebra $A$ has  infinite homological dimension, then we are forced to deal
with the derived category of right  bounded  complexes (in order
to define  the left derived functor of the tensor product). In case $A$ has  
finite homological dimension we can suppose that all complexes above are {\it bounded}
from  both sides.
\end{remark}

\begin{theorem}
The functor 
$$
D^{-}(A-\mod) \stackrel{\mathbf F}\lar {\rm TC}_{A}
$$
$
\kP_{\bullet} \lar
(\tilde{A}\tens_{A} {\mathcal P}_{\bullet} ,
A/I\tens_{A}{\mathcal P}_{\bullet},
i : A/I\tens_{A} {\mathcal P}_{\bullet}   \lar \tilde{A}/I\tens_{A}{\mathcal P}_{\bullet}) 
$
has the following properties:

\begin{enumerate}
\item ${\mathbf F}$ is dense {\rm (}i.e., every triple 
$(\tilde{\kP}_{\bullet}, \kM_{\bullet}, i)$ is 
isomorphic to some
${\mathbf F}(\kP_{\bullet})${\rm )}. 
\item $\FF(\kP_{\bullet}) \cong \FF(\kQ_{\bullet}) \Longleftrightarrow 
\kP_{\bullet} \cong \kQ_{\bullet}$.
\item ${\mathbf F}(\kP_{\bullet})$ is indecomposable if and only if so  is $\kP_{\bullet}$
(note that this property is an easy formal consequence of the previous two properties).
\item ${\mathbf F}$ is full.
\end{enumerate} 
\end{theorem}

\begin{remark} ${\mathbf F}$ is not faithful. So it is not an equivalence
of categories. A functor  ${\mathbf F}$ satisfying the properties
1-4 is called  \emph{detecting functor} (see \cite{Baues}). 
\end{remark}

\noindent
{\it Proof.}
The main point to be clarified is: having a triple 
$\kT = (\tilde{\kP}_{\bullet}, \kM_{\bullet}, i)$ how can we reconstruct $\kP_{\bullet}$?
The exact sequence 
$$0\lar I\tilde{\kP}_{\bullet} \lar \tilde{\kP}_{\bullet} \lar
\tilde{A}/I\tens_{\tilde{A}} \tilde{\kP}_{\bullet} \lar 0$$ of complexes
in $A-\mod$ gives a distinguished triangle 
$$
I\tilde{\kP}_{\bullet} \lar \tilde{\kP}_{\bullet} \lar
\tilde{A}/I\tens_{\tilde{A}} \tilde{\kP}_{\bullet}
\lar I\tilde{\kP}_{\bullet}[-1]
$$ 
in $D^{-}(A-\mod)$.
The properties of triangulated categories imply that there is a
 morphism of triangles

\begin{tabular}{p{1.1cm}c}
&
\xymatrix
{
I\tilde{\mathcal P}_{\bullet} \ar[r] &
\tilde{\mathcal P}_{\bullet} \ar[r] &
\tilde{A}/I\tens_{\tilde{A}} \tilde{\mathcal P}_{\bullet}  \ar[r]&
I\tilde{\mathcal P}_{\bullet}[-1] \\
I\tilde{\mathcal P}_{\bullet} \ar[r] \ar[u]^{id} &
{\mathcal P}_{\bullet} \ar[r] \ar[u]_{\Phi} &
{\mathcal M}_{\bullet} \ar[r] \ar[u]_{i} &
I\tilde{\mathcal P}_{\bullet}[-1],\ar[u]^{id}\\
}
\end{tabular}

\noindent
where $\kP_{\bullet} = cone(\kM_{\bullet} \lar I \tilde{\kP_{\bullet}}[-1])[1]$. Set ${\bf G}(\kT) = \kP_{\bullet}.$ Taking a cone
is not a functorial operation. It gives an intuitive explanation why
the functor ${\bf F}$ is not an equivalence.
The properties of triangulated categories  immediately imply that
the constructed map (not a functor!)
$$
{\bf G} : {\rm Ob}({\rm TC}_{A}) \lar {\rm Ob}\bigl(D^{-}(A-\mod)\bigr)
$$
sends isomorphic objects into isomorphic ones and 
${\bf GF}(\kP_{\bullet}) \cong \kP_{\bullet}$. 
Now we have to show that
${\bf FG}(\tilde{\kP}_{\bullet}, \kM_{\bullet}, i) \cong (\tilde{\kP}_{\bullet}, \kM_{\bullet}, i)$.

\begin{lemma}
In the above notations,  let
$\tilde{P}$ be a projective $\tilde{A}$-module, $M$ be an $A/I$-module,
$i: M\lar \tilde{P}/I\tilde{P}$ be an $A/I$-module monomorphism such that
the induced map $\tilde{i}: \tilde{A}/I \tens_{A/I} M \lar \tilde{P}/I\tilde{P}$
is an isomorphism. Consider the pull-back diagram 

\begin{tabular}{p{1.8cm}c}
&
\xymatrix
{0 \ar[r] & I \tilde{P} \ar[r] \ar[d]^{id} & P \ar[r] \ar[d] &
 M \ar[r] \ar[d]^{i} & 0 \\
0 \ar[r] & I \tilde{P} \ar[r]  & \tilde{P} \ar[r]^{\pi}  &
 \tilde{P}/ I \tilde{P} \ar[r] & 0.
}
\end{tabular}

\noindent
Then $P$ is a projective $A$-module and $\tilde{A}\tens_{A} P \lar \tilde{P}$
is an isomorphism.
\end{lemma}

Consider the image $\bar{I}$ of the ideal $I$ in $\tilde{A}/r$. Since 
$\tilde{A}/r$
is semi-simple, we can find an ideal $\bar{J}$ in $\tilde{A}/r$ such that 
$\bar{I} + \bar{J} = \tilde{A}/r$, $\bar{I}\cap \bar{J} = 0$. By the Chinese
remainder theorem we have $\tilde{A}/r = \tilde{A}/I \times \tilde{A}/J$.

Let $\overline{\tilde{P}} = \tilde{P}/r\tilde{P}$. Then $\overline{\tilde{P}} =
\overline{\tilde{P}}_1\oplus \overline{\tilde{P}}_2$, where $\overline{\tilde{P}}_1$ is
an $\tilde{A}/I$-module and $\overline{\tilde{P}}_2$ an $\tilde{A}/J$-module.
But then $\tilde{P}$ also decomposes into a direct sum:
$\tilde{P} = \tilde{P}_{1} \oplus \tilde{P}_{2}$, where
$\overline{\tilde{P}}_{i} = \tilde{P}_{i}/
r \tilde{P}_{i}, i = 1,2$ (we use the fact
that there is a bijection between projective and semi-simple modules:
$\tilde{P} \longleftrightarrow \tilde{P}/r\tilde{P}$). 

\noindent
Then we have:
$$
I\tilde{P}_{1} = r\tilde{P}_{1}, 
I\tilde{P}_{2} = \tilde{P}_{2}.
$$
Indeed, $\tilde{P}_{1}/r \tilde{P}_{1}$ is an $A/I$-module, so 
$I \tilde{P}_{1} \subseteq r \tilde{P}_{1}$. But $r\subseteq I$, hence
$r \tilde{P}_{1} \subseteq I \tilde{P}_{1}$. So, 
$I\overline{\tilde{P}}_{1} = r \overline{\tilde{P}}_{1}$. 
Analogously, $J \tilde{P}_{2} \subseteq r \tilde{P}_{2}$. But $I+J = \tilde{A}$,
so $$\tilde{P}_{2}= I\tilde{P}_{2} +  J \tilde{P}_{2} \subseteq 
I\tilde{P}_{2} +  r \tilde{P}_{2} \subseteq \tilde{P}_{2}.$$

\noindent
Hence, by Nakayama's Lemma $I\tilde{P}_{2} = \tilde{P}_{2}$.

Our diagram has now the form:

\begin{tabular}{p{0.7cm}c}
&
\xymatrix
{0 \ar[r] & \tilde{P}_{2} \oplus r \tilde{P}_{1} \ar[r] \ar[d]^{id} & P \ar[r] \ar[d] &
 M \ar[r] \ar[d]^{i} & 0 \\
0 \ar[r] & \tilde{P}_{2} \oplus r \tilde{P}_{1} \ar[r]  & \tilde{P}_{1} \oplus 
\tilde{P}_{2} \ar[r]^{\pi}  &
 \tilde{P}_{1}/ r \tilde{P}_{1} \ar[r] & 0.
}
\end{tabular}

\noindent
Since $P \lar \tilde{P}_{1} \oplus \tilde{P}_{2}$ is a monomorphism,
$\tilde{P}_{2}$ is a direct summand of $P$. Moreover, $\tilde{P}_{2}$ is
a projective $A$-module. 
Indeed, let $\tilde{Q}$ be any projective $\tilde{A}$-module satisfying
$I\tilde{Q} = \tilde{Q}$. Without loss of generality  suppose that $\tilde{Q}$
 is a direct summand of $\tilde{A}$. Then 
$$\tilde{Q} = I\tilde{Q} \subseteq I\tilde{A} \subseteq A \subseteq \tilde{A}.$$
But if the embedding $\tilde{Q} \lar \tilde{A}$ splits, then 
$\tilde{Q}\lar A$ splits too.
Hence,  $\tilde{Q}$ is a projective $A$-module.

Note that the canonical map $\tilde{A}\tens_{A} \tilde{P}_{2} \lar 
\tilde{P}_{2}$ is an isomorphism.  Indeed,
$$\tilde{A}\tens_{A} \tilde{P}_{2} = \tilde{A}\tens_{A} I\tilde{P}_{2} =
\tilde{A} I \tens_{A} \tilde{P}_{2} = I \tens_{A} \tilde{P}_{2}.$$
But $\tilde{P}_{2}$ is a flat $A$-module, hence
 $$I \tens_{A} \tilde{P}_{2} = 
I\tilde{P}_{2} = \tilde{P}_{2}.$$

\noindent
So we get

\begin{tabular}{p{1.5cm}c}
&
\xymatrix
{0 \ar[r] & r \tilde{P}_{1} \ar[r] \ar[d]^{id} & P_{1} \ar[r] \ar[d] &
 M \ar[r] \ar[d]^{i} & 0 \\
0 \ar[r] & r \tilde{P}_{1} \ar[r] & \tilde{P}_{1}  \ar[r]^{\pi}  &
 \tilde{P}_{1}/ r \tilde{P}_{1} \ar[r] & 0.
}
\end{tabular}

We know that  
$\tilde{i}: \tilde{A}/I\tens_{A/I} M \lar \tilde{P}_{1}/ I \tilde{P}_{1}$
is an isomorphism. But then 
$\tilde{A}/r\tens_{A/r} M \lar \tilde{P}_{1}/ r \tilde{P}_{1}$ is 
an isomorphism, too. Indeed  $IM = 0,$ since $M$ is a submodule of $\tilde{P}_{1}/ I \tilde{P}_{1}.$ But $I+J = A$, hence $JM = M$ and
$\tilde{A}/J\tens_{A/r} M = 0$. Therefore 
$$ \tilde{A}/I\tens_{A/I} M = \tilde{A}/I\tens_{A/r} M \cong 
(\tilde{A}/I \oplus \tilde{A}/J)\tens_{A/r} M \cong \tilde{A}/r\tens_{A/r} M.$$

Now we have  to show that
$P_{1}$ is projective and 
$\tilde{A}\tens_{A} P \lar \tilde{P}_{1}$ is an 
isomorphism. 
Let $P(M)$ be a projective covering of $M$.

\begin{tabular}{p{1.6cm}c}
&
\xymatrix
{
& & & P(M) \ar[d] \ar[dl]_{\psi} & \\
0 \ar[r] &r\tilde P_{1} \ar[r] \ar[d]^{id} & P_{1} \ar[r] \ar[d] &
M \ar[r] \ar[d]^{i} & 0 \\
0 \ar[r] &r\tilde{P}_{1} \ar[r]  & \tilde{P}_{1} \ar[r]  &
\tilde{P_{1}}/r\tilde{P}_{1} \ar[r] & 0. \\
}
\end{tabular}

Apply the functor $A/r\tens_{A} $ to the first row of this diagram. We get:
$\bar{\psi}: P(M)/rP(M) \lar P/rP$ is an isomorphism. Hence by Nakayama's lemma
 $\psi$ is an epimorphism.
Consider the composition map $P(M) \lar \tilde{P}_{1}$. The induced map 
$\tilde{A}\tens_{A} P(M) \lar \tilde{P}_{1}$ is an isomorphism modulo $r$. 
Since both modules are projective, it is indeed an isomorphism. 
We get: $P(M) \lar \tilde{A}\tens P(M) \lar \tilde{P}_{1}$ is a monomorphism.
But then $\psi: P(M) \lar P_1$ is a monomorphism too. So it is an isomorphism.
And we have shown also that $\tilde{A}\tens_{A} P \lar \tilde{P}_{1}$ is an 
isomorphism.

We finish now the proof of the theorem. Let $(\tilde{\kP}_{\bullet}, \kM_{\bullet}, i)$ be
a triple. Without loss of generality, suppose that $\tilde{\kP}_{\bullet}$ is a minimal
complex and $\kM_{\bullet}$ a complex with zero differentials. Then
$\tilde{\kP}_{\bullet}/ r \tilde{\kP}_{\bullet}$ is a  complex with zero differentials too and 
the map
$i: \kM_{\bullet} \lar \tilde{\kP}_{\bullet}/ I \tilde{\kP}_{\bullet}$ has the  property that
$\tilde{i}: \tilde{A}/I_{\tens{A/I}}  \kM_{\bullet} \lar \tilde{\kP}_{\bullet}/ I \tilde{\kP}_{\bullet}$ is an isomorphism of {\it complexes}.  Consider now
the pull-back diagram in the abelian category of complexes of $A$-modules.

\begin{tabular}{p{1.6cm}c}
&
\xymatrix
{0 \ar[r] & I\tilde{\kP}_{\bullet} \ar[r] \ar[d]^{id} & \kP_{\bullet} \ar[r]^{\Psi} 
\ar[d]^{\Phi} &
 \kM_{\bullet} \ar[r] \ar[d]^{i} & 0 \\
0 \ar[r] & I\tilde{\kP}_{\bullet} \ar[r]  & \tilde{\kP}_{\bullet} \ar[r]^{\pi}  &
 \tilde{\kP}_{\bullet}/  I\tilde{\kP}_{\bullet} \ar[r] & 0.
}
\end{tabular}

\noindent
From Lemma 2.6 follows that: 
\begin{enumerate}
\item $\kP_{\bullet}$ is a complex of projective $A$-modules;
\item $(\Phi\tens id, \Psi\tens id): 
(\tilde{A}\tens_{A}\kP_{\bullet}, 
A/I\tens_{A}\kP_{\bullet},
A/I\tens_{A} \kP_{\bullet} \lar 
\tilde{A}\tens_{A}\kP_{\bullet}) \lar 
(\tilde{\kP}_{\bullet}, \kM_{\bullet}, i)$
is an isomorphism in the category of triples. 
\end{enumerate}

\noindent
It remains to show that {\bf F}
is full.

Let $(\Phi, \varphi) : (\tilde{\kP_{\bullet 1}}, \kM_{\bullet 1}, i_1) \lar 
(\tilde{\kP_{\bullet 2}}, {\kM_{\bullet 2}}, i_2)$ be a morphism in 
$\TC_{A}$, where $\kM_{\bullet 1}$ and $\kM_{\bullet 2}$ are complexes with zero differentials.
Since we are dealing with complexes of projective objects, $\Phi$ and $\varphi$
can be represented by morphisms of complexes. Let us moreover suppose 
$\tilde{\kP}_{\bullet 1}$ and $\tilde{\kP}_{\bullet 2}$ to be minimal.
Then $\tilde{A/I}\tens_{\tilde{A}} \tilde{\kP}_{\bullet 1}, i=1,2$ are complexes
with zero differentials, too.

\begin{tabular}{p{2.7cm}c}
&
\xymatrix
{ \mathcal M_{\bullet 1} \ar[r]^{\varphi} \ar[d]^{i_1} & 
\mathcal M_{\bullet 2}
\ar[d]^{i_2} \\
  \tilde{\mathcal P}_{\bullet 1}/ I\tilde{\kP}_{\bullet 1} \ar[r]^{\bar{\Phi}} &
  \tilde{\kP}_{\bullet 2}/ I\tilde{\kP}_{\bullet 2} \\
}
\end{tabular}

\noindent
is commutative 
{\it in the category of complexes}. The properties of pull-back imply the 
existence of a  morphism of complexes 
$ \kP_{\bullet 1} \lar \kP_{\bullet 2} $ such that  
$$
\begin{tabular}{p{0.2cm}c}
&
\xymatrix@!0
{ & \kP_{\bullet 2} \ar[rr]\ar'[d][dd] & & \kM_{\bullet 2} \ar[dd] \\
 \kP_{\bullet 1} \ar[ur] \ar[rr] \ar[dd]& & \kM_{\bullet 1} \ar[ur]
\ar[dd] \\
& \tilde{\kP}_{\bullet 2} \ar[rr] & 
& \tilde{\kP}_{\bullet 2}/ I\tilde{\kP}_{\bullet 2} \\
\tilde{\kP}_{\bullet 1} \ar[rr] \ar[ur] & 
& \tilde{\kP}_{\bullet 1}/ I\tilde{\kP}_{\bullet 1}\ar[ur]
}
\end{tabular} 
$$
is commutative.  Hence it gives a lift of a morphism $(\Phi, \varphi)$ we are 
looking for. So the functor ${\bf F}$ is full, which accomplishes the proof of the Theorem.

\section{The case of 
$D^{-}(\kk\langle\langle x,y \rangle\rangle /(x^2,y^2))$}

Consider the embedding of completed path algebras $A \lar \tilde{A}$:

\begin{figure}[ht]
\hspace{3.3cm}
\includegraphics[height=1.5cm,width=5cm,angle=0]{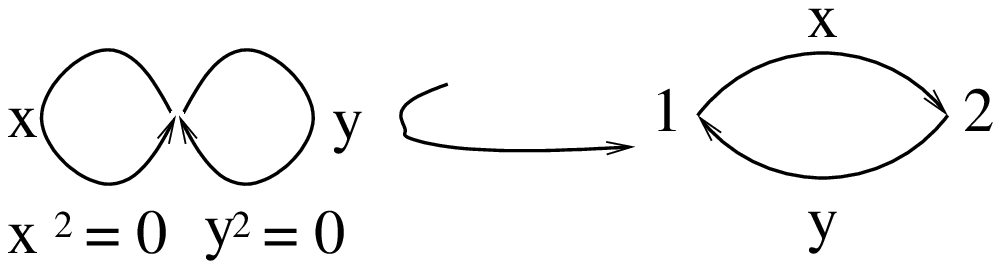}
\end{figure}

Take  $I = (x,y)$, then 
$A/I = \kk, \tilde{A}/I = \kk\times \kk$ and $A/I \lar \tilde{A}/I $ is 
just the diagonal map. 

As we have seen in the previous section, a complex $\kP_{\bullet}$ of the
derived category $D^{-}(A-\mod)$ is defined by some triple
$(\tilde{\kP}_{\bullet},\kM_{\bullet},i)$. Since $A/I-\mod$ can be identified with
the category of $\kk$-vector spaces, the map 
$i: \kM_{\bullet} \lar \tilde{\kP}_{\bullet}/I\tilde{\kP}_{\bullet}$ is given by a collection
of linear maps
$H_{k}(i): H_{k}(\kM_{\bullet}) \lar H_{k}(\tilde{\kP}_{\bullet}/I\tilde{\kP}_{\bullet})$.
The map $H_{k}(i)$ is a $\kk$-linear map of a $\kk$-module into a 
$\kk\times \kk$-module.  Hence it is given by two matrices $H_{k}(i|1)$ and 
$H_{k}(i|2)$. From the non-degeneracy condition of the category of
 triples it follows that both of these matrices are invertible.

The algebra $\tilde{A}$ has  homological dimension $1$. Moreover,  it is an
order. By a  theorem of Dold (see \cite{Dold}), an indecomposable complex
from $D^{-}(\tilde{A}-\mod)$ is isomorphic to
$$
\dots \lar 0 \lar \underbrace{M}_{i} \lar 0 \lar \dots, 
$$ 
where $M$ is an indecomposable $\tilde{A}$-module. 
But $\tilde{A}$ is a hereditary order over $\kk[[t]]$. Hence (see \cite{uniserial}) an indecomposable finitely-generated  $\tilde A$-module is isomorphic in the derived 
category to a shift of $P_{i}$, $i =1,2$ or to 
$$
P_{i} \stackrel{\varphi}\lar P_{j}  \qquad 
(i,j = 1,2),
$$ 
where $P_{i} = \tilde{A}e_{i}$ ($e_i$ is the idempotent corresponding to the point $i$ of the graph), $\varphi$ is a morphism given by a path  going from  $j$ to $i$ and ${\rm Im}(\varphi) \subseteq rP_{j}$, thus
 $$\tilde{A}/I \tens_{\tilde{A}} ({P}_{i} \stackrel{\varphi}\lar P_{j}) = 
\kk_{i} \stackrel{0}\lar \kk_{j}.$$
Let 
$$\tilde{\kP}_{\bullet} = \bigoplus \tilde{\kP}_{\bullet i}^{n_{i}}$$
be a  decomposition of $\tilde{\kP}_{\bullet}$ into a direct sum of indecomposables. This decomposition
implies a division of matrices $H_{k}(i|1)$ and $H_{k}(i|2)$ into horizontal
stripes.

The next question is: which transformations can we perform with the matrices
$H_{k}(i|1)$ and $H_{k}(i|2)$? 

We can do simultaneously any elementary transformation of columns of 
$H_{k}(i|1)$ and $H_{k}(i|2)$ (they correspond to the automorphisms of
$\kM_\bullet$). From the definition of the category of triples
it follows that row transformations are induced by  morphisms in
$D^{-}(\tilde{A}-\mod)$.

Let us now describe the morphisms between indecomposable complexes
from $D^{-}(\tilde{A}-\mod)$, which are non-zero after applying 
$\tilde{A}/I$. Due to  \cite{uniserial} they are just

\begin{tabular}{p{2.2cm}ccc}
&
\xymatrix
{
P_{2} \ar[r]^{\cdot x} & P_{1}\\
P_{1} \ar[r]^{\cdot yx} \ar[u]^{\cdot y} & P_{1} \ar[u]\\
P_{2} \ar[r]^{\cdot xyx} \ar[u]^{\cdot x} & P_{1} \ar[u]\\
P_{1} \ar[r]^{\cdot yxyx} \ar[u]^{\cdot y} & P_{1} \ar[u]\\
\dots \ar[u] & \dots \ar[u]\\
}
&
\qquad
&
\xymatrix
{
P_{1} \ar[r]^{\cdot y} \ar[d] & P_{2} \ar[d]^{\cdot x}\\
P_{1} \ar[r]^{\cdot yx} \ar[d] & P_{1} \ar[d]^{\cdot y}\\
P_{1} \ar[r]^{\cdot yxy} \ar[d] & P_{2} \ar[d]^{\cdot x}\\
P_{1} \ar[r]^{\cdot yxyx} \ar[d] & P_{1} \ar[d]\\
\dots  & \dots \\
}
\end{tabular}

\begin{tabular}{p{2.2cm}ccc}
&
\xymatrix
{
P_{1} \ar[r]^{\cdot y} & P_{2}\\
P_{2} \ar[r]^{\cdot xy} \ar[u]^{\cdot x} & P_{2} \ar[u]\\
P_{1} \ar[r]^{\cdot yxy} \ar[u]^{\cdot y} & P_{2} \ar[u]\\
P_{2} \ar[r]^{\cdot xyxy} \ar[u]^{\cdot x} & P_{2} \ar[u]\\
\dots \ar[u] & \dots \ar[u]\\
}
&
\qquad
&
\xymatrix
{
P_{2} \ar[r]^{\cdot x} \ar[d] & P_{1} \ar[d]^{\cdot y}\\
P_{2} \ar[r]^{\cdot xy} \ar[d] & P_{2} \ar[d]^{\cdot x}\\
P_{2} \ar[r]^{\cdot xyx} \ar[d] & P_{1} \ar[d]^{\cdot y}\\
P_{2} \ar[r]^{\cdot xyxy} \ar[d] & P_{2} \ar[d]\\
\dots  & \dots \\
}
\end{tabular}

\begin{tabular}{p{3.8cm}ccc}
&
\xymatrix
{
0 \ar[r] \ar[d] & P_{i} \ar[d]\\
P_{j} \ar[r]^{\varphi}  & P_{i} \\
}
\end{tabular}

\begin{tabular}{p{3.8cm}ccc}
&
\xymatrix
{
P_{j} \ar[r]^{\varphi} \ar[d] & P_{i}  \ar[d]\\
P_{j} \ar[r] & 0\\
}
\end{tabular}

Moreover, we always have a morphism

\begin{tabular}{p{3cm}ccc}
&
\xymatrix
{
0 \ar[r] \ar[d] &P_{i} \ar[d]^{\lambda\cdot id} \ar[r] & P_{j} \ar[d]\\
P_{k} \ar[r]  & P_{i} \ar[r] & 0.\\
}  
\end{tabular}

Now note that we have the following  cases:
\begin{enumerate}
\item  A morphism 

\begin{tabular}{p{2.7cm}ccc}
&
\xymatrix
{
0 \ar[r] \ar[d] &P_{i} \ar[d]^{\lambda\cdot id} \ar[r] & P_{j} \ar[d]\\
P_{k} \ar[r]  & P_{i} \ar[r] & 0\\
}  
\end{tabular}

\noindent
induces

\begin{tabular}{p{2.7cm}ccc}
&
\xymatrix
{
0 \ar[r] \ar[d] & \kk_{i} \ar[d]^{\lambda} \ar[r]^0 & \kk_{j} \ar[d]\\
\kk_{k} \ar[r]^0  & \kk_{i} \ar[r] & 0.\\
}  
\end{tabular}

\item A morphism

\begin{tabular}{p{3.4cm}cc}
&
\xymatrix
{
 P_{j}  \ar[r] \ar[d]^{\varphi} & P_{i} \ar[d]^{\lambda\cdot id}\\
P_{k} \ar[r]  & P_{i} \\
}  
\end{tabular}

\noindent
where $\varphi \in {\rm rad}(\tilde{A})$,
induces

\begin{tabular}{p{3.4cm}cc}
&
\xymatrix
{
\kk_{j} \ar[d]^{0} \ar[r] & \kk_{j} \ar[d]^{\lambda}\\
\kk_{k} \ar[r]^0 &  \kk_{i}.\\
}  
\end{tabular}

\item Analogously we have that

\begin{tabular}{p{3.4cm}cc}
&
\xymatrix
{
 P_{i}  \ar[d]^{\lambda\cdot id} \ar[r] & P_{j} \ar[d]^{\varphi}\\
P_{i} \ar[r]  & P_{k} \\
}  
\end{tabular}

\noindent
where $\varphi \in {\rm rad}(\tilde{A})$,
induces

\begin{tabular}{p{3.5cm}cc}
&
\xymatrix
{
\kk_{i} \ar[d]^{\lambda} \ar[r]^0 & \kk_{j} \ar[d]^{0}\\
\kk_{i} \ar[r]^0 &  \kk_{k}.\\
}  
\end{tabular}

\item In the same way

\begin{tabular}{p{3.5cm}ccc}
&
\xymatrix
{
0 \ar[r] \ar[d] & P_{i} \ar[d]^{\lambda\cdot id}\\
P_{j} \ar[r]  & P_{i} \\
}
\end{tabular}

\noindent
induces 

\begin{tabular}{p{3.5cm}ccc}
&
\xymatrix
{
0 \ar[r] \ar[d] & \kk_{i} \ar[d]^{\lambda}\\
\kk_{j} \ar[r]^{0}  & \kk_{i} \\
}
\end{tabular}

\noindent and

\begin{tabular}{p{3.5cm}ccc}
&
\xymatrix
{
P_{j} \ar[r]^{\varphi} \ar[d]^{\lambda\cdot id} & P_{i} \ar[d]\\
P_{j} \ar[r] & 0\\
}
\end{tabular}

\noindent induces 

\begin{tabular}{p{3.5cm}ccc}
&
\xymatrix
{
\kk_{j} \ar[r]^{0} \ar[d]^{\lambda} & \kk_{i} \ar[d]\\
\kk_{j} \ar[r] & 0.\\
}
\end{tabular}

\item And finally an endomorphism

\begin{tabular}{p{3.5cm}cc}
&
\xymatrix
{
 P_{j}  \ar[r]^{\varphi} \ar[d]^{\lambda\cdot id} & P_{i} \ar[d]^{\lambda\cdot id}\\
 P_{j}  \ar[r]^\varphi  & P_{i} \\
}  
\end{tabular}

\noindent
induces 

\begin{tabular}{p{3.5cm}cc}
&
\xymatrix
{
\kk_{i} \ar[d]^{\lambda} \ar[r]^0 & \kk_{j} \ar[d]^{\lambda}\\
\kk_{i}  \ar[r]^0 & \kk_{j},\\
}  
\end{tabular}

\noindent
and the same for $P_{i}, i=1,2$. 
\end{enumerate}

From what has been said we observe that the matrix problem describing the derived category $D^{-}(A-\mod)$ is given by the following partially ordered set
(bunch of chains, see \cite{mp} or Appendix B \cite{vb}).

\begin{figure}[ht]
\hspace{2.8cm}
\includegraphics[height=8.5cm,width=6cm,angle=0]{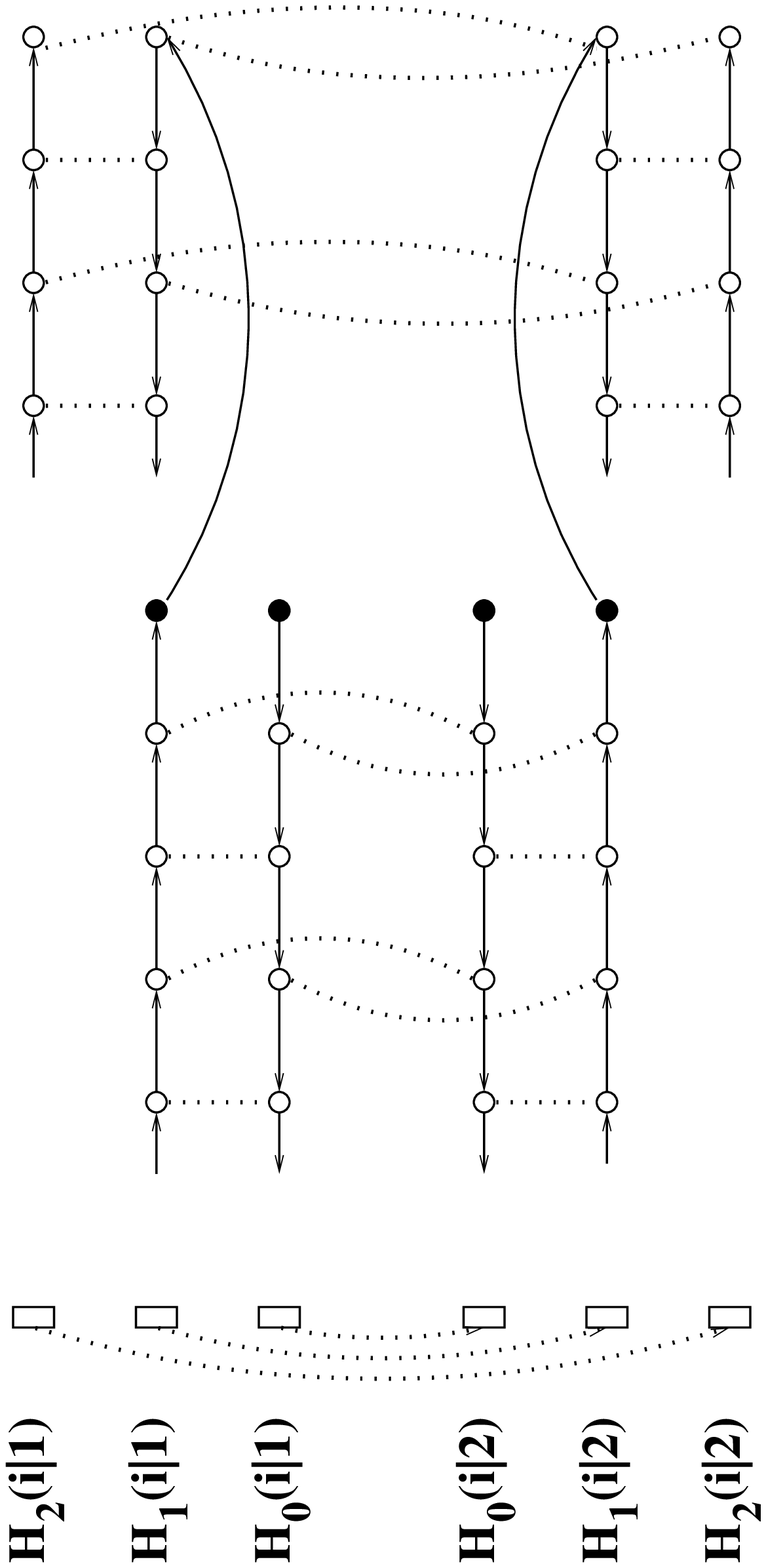}
\end{figure}

In this picture we assume that complexes are shifted in such a way that 
all $H_k(i) = 0$ for $i < 0$. 
 Small circles correspond to the horizontal stripes, small rectangles correspond to the vertical stripes, dotted lines between
circles show
the related stripes (i.e. those which come from the same object of the derived category),  vertical arrows 
describe the possible transformations between different horizontal stripes:

Explicitly saying, we can do the following  transformations with our matrices
$H_{\bullet}(i)$:

\begin{enumerate}
\item We can do any simultaneous elementary transformations of the columns of 
the matrices
 $H_{k}(i|1)$ and $H_{k}(i|2), k\in \mathbb Z$.
\item We can do any simultaneous transformations of rows inside conjugated blocks.
\item We can add a scalar multiple of any row from a block 
with lower weight to any row of a block of a higher weight (inside 
 the big matrix, of course). These transformations can be proceeded 
independently
inside 
 $H_{k}(i|1)$ and $H_{k}(i|2)$, $k\in \mathbb Z$.
\end{enumerate}

This matrix problem belongs to the well-known representations of bunches of 
chains (see \cite{mp}, \cite{Roiter3}, \cite{Bond} and Appendix). From here 
we conclude that
there are 
 three types of indecomposable objects: bands, finite strings (both correspond to complexes of finite projective dimension) and infinite strings 
(which correspond to complexes of infinite projective  dimension). 
In the  section 6 we shall explain,
how the combinatoric of band and 
string representations can be used to write down explicit projective resolutions of complexes.

\section{Gelfand quiver}
In this section we shall see that our technique allows us to describe 
the derived category of representations of the completed path algebra
of the quiver

\begin{figure}[ht]
\hspace{4cm}
\includegraphics[height=1.5cm,width=5cm,angle=0]{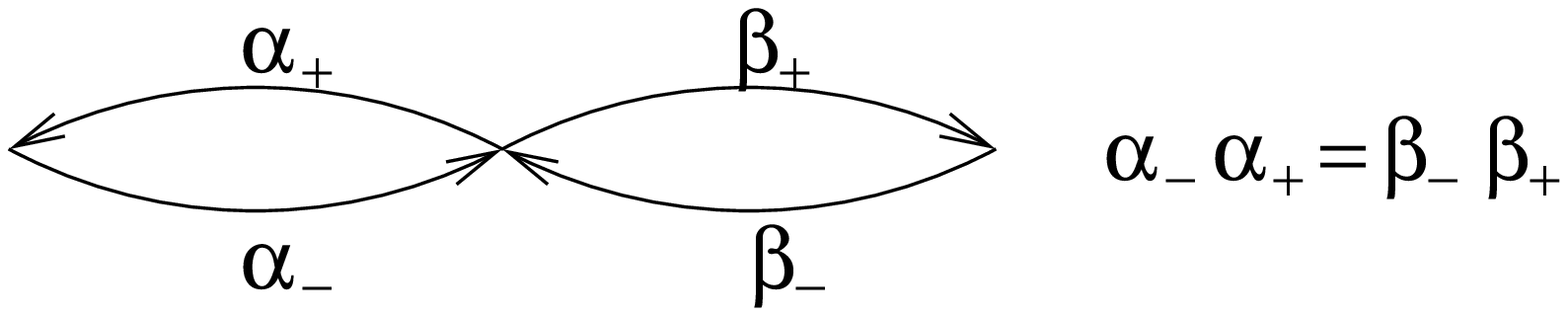}
\end{figure}

The classification of indecomposable representations of this quiver
can be  reduced to representations of bunches of {\it semi-chains}, see
\cite{mp}. It is
not surprising that the description of the derived category is reduced to
the problem of the same type. Consider the embedding given in Example 1.3.
In this case we have: $A/I = \kk\times \kk$, $\tilde{A}/I = M_{2}(\kk)$ and
$A/I \lar \tilde{A}/I$ the diagonal mapping. Now we have to answer
the following:
\begin{question}
Let $M$ be a $\kk\times \kk$-module, $M'$ be a $M_{2}(\kk)$-module, 
$\varphi: M \lar M'$ a map of $\kk\times \kk$-modules ($M'$ is supplied with 
$\kk\times \kk$-module structure using the diagonal embedding). The map
of $\kk\times \kk$-modules is given by two matrices $\varphi(1)$ and $\varphi(2)$. 
Which conditions
should satisfy $\varphi(1)$ and $\varphi(2)$ in order 
$\tilde{i}:  M_{2}(\kk)\tens_{\kk \times \kk} M \lar M'$ to be an isomorphism?
\end{question}
Let 
$M = \langle v_{1}, v_{2},\dots,v_{m}; w_{1}, w_{2},\dots,w_{n} \rangle =
\kk(1)^{m} \oplus \kk(2)^{n}.$
There is only one indecomposable $M_{2}(\kk)$-module:  $\kk^{2}$. So, 
$$
 M' = \langle u_{1}',u_{1}'';  u_{2}',u_{2}''; \dots,  u_{N}',u_{N}'' \rangle =
(\kk^{2})^{N},
$$
where the action of matrix units are:
$$
e_{11} u_{i}' = u_{i}',\quad e_{21} u_{i}' = u_{i}'', \quad
e_{12} u_{i}' = 0, \quad e_{22} u_{i}' = 0,
$$
and, analogously,
$$
e_{11} u_{i}'' = 0, \quad e_{21} u_{i}''= 0, \quad
e_{12} u_{i}'' = u_{i}', \quad e_{22} u_{i}'' =  u_{i}''.
$$
Let
$$
\varphi(v_{i}) = \sum_{j=1}^{N} \alpha_{ji} u_{j}' + 
\sum_{j=1}^{N} \alpha'_{ji} u_{j}''.
$$ 
Since $\varphi$ is a $\kk\times \kk$-module homomorphism, 
$0 = \varphi(e_{22}v_{i}) = e_{22} \varphi(v_{i})$. So 
 all $\alpha'_{ji} = 0$. 
Analogously,
$$
\varphi(w_{i}) = \sum_{j=1}^{N} \beta_{ji} u_{j}''.
$$

On the other hand, any $M_{2}(\kk)$ homomorphism 
$
\psi : (\kk^{2})^{n} \lar (\kk^{2})^{m}
$
is given by an $m\times n$ matrix $(\alpha_{ij})$ with the entries from $\kk$
(see \cite{DrozdKirich}, theorem 1.7.5). Namely, if 
$(e'_{1}, e''_{1},e'_{2}, e''_{2}, \dots,
e'_{n},e''_{n})$ and  $(f'_{1}, f''_{1},f'_{2}, f''_{2}, \dots,
f'_{m},f''_{m})$ are canonical bases of $(\kk^{2})^{n}$ and  $(\kk^{2})^{m}$
then 
$$
\psi(e'_{j}) = \sum_{j=1}^{m} \alpha_{ij} f'_{i},
$$
$$
\psi(e''_{j}) = \sum_{j=1}^{m} \alpha_{ij} f''_{i}.
$$
Consider now a $M_{2}(\kk)$-module $M_{2}(\kk)\tens_{\kk\times \kk} M$. It is
generated by
$$
e_{11}\tens v_{1}, e_{21}\tens v_{1};
e_{11}\tens v_{2}, e_{21}\tens v_{2};\dots; 
e_{11}\tens v_{m}, e_{21}\tens v_{m};
$$
$$
e_{12}\tens w_{1}, e_{22}\tens w_{1};
e_{12}\tens w_{2}, e_{22}\tens w_{2}; \dots;
e_{12}\tens w_{n}, e_{22}\tens w_{n}.
$$
Since $\tilde{\varphi}(e\tens v) = e\varphi(v)$, it is easy to see that
$\varphi$ is given by $N\times (n+m)$-matrix
 $\bigl(\varphi(1)| \varphi(2)\bigr)$.
So, $\tilde{\varphi}$ is an isomorphism if 
$\bigl(\varphi(1)| \varphi(2)\bigr)$ is
square and invertible. 

Now let us return to the Gelfand quiver. The morphisms in
$D^{-}(\tilde{A}-\mod)$ were discussed in the previous subsection. 
We are able to write Bondarenko's partially ordered set:

\begin{figure}[ht]
\hspace{3.9cm}
\includegraphics[height=6.5cm,width=4cm,angle=0]{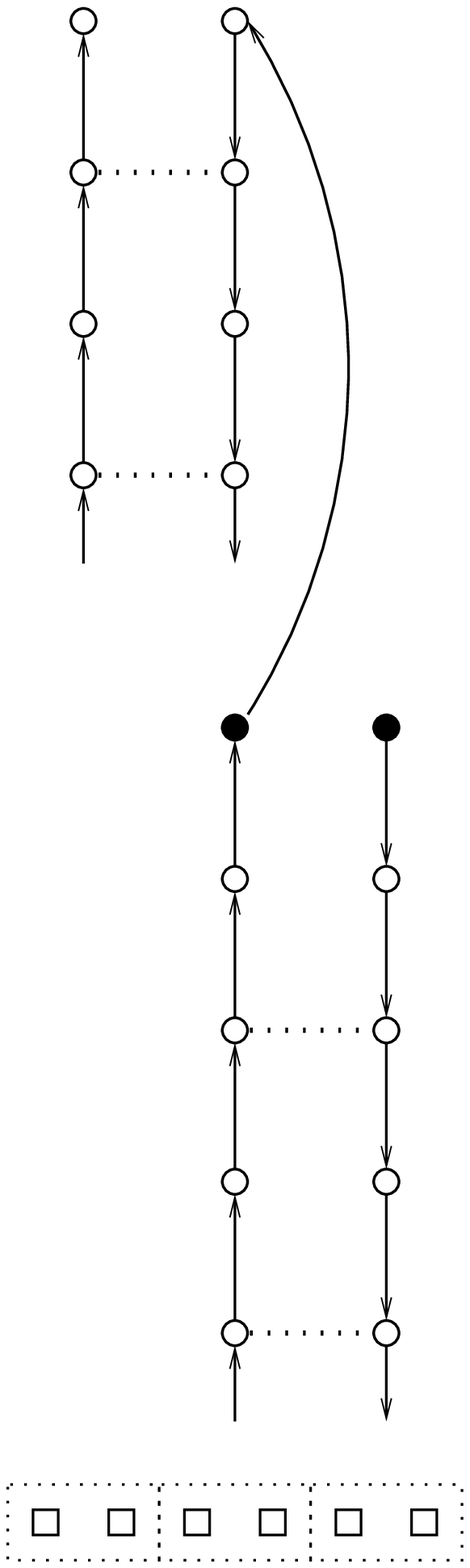}
\end{figure} 

This picture shows the division of matrices $H_{k}(i), k\in \mathbb Z$, into horizontal
and vertical stripes. Each of these matrices is divided into two vertical
blocks $H_{k}(i|1)$ and $H_{k}(i|2)$ (which correspond to the fact that we have an embedding 
$\kk\times \kk \lar M_{2}(\kk) )$ and horizontal blocks  that correspond to
indecomposables of $D^{-}(\tilde{A}-\mod)$. In the same way as in the previous 
section we have an  ordering on the horizontal stripes.

We can perform the following transformations with matrices
$H_{\bullet}(i)$:

\begin{enumerate}
\item We can do independently elementary transformations of  columns 
 of   $H_{k}(i|1)$ and $H_{k}(i|2)$.
\item We can do any simultaneous transformations of rows inside conjugated blocks.
\item We can add a scalar multiple of any row from a block 
with lower weight to any row of a block of a higher weight.
\end{enumerate}

This problem belongs to the class of representations of bunches of
 {\it semi-chains}. The description of indecomposable objects was 
obtained in  \cite{mp} and \cite{CB} and later elaborated in \cite{Deng}. 
Since we get in this case infinitely many matrices, certain modifications 
should be done, see \cite{Bond} and Appendix. 
Namely, there are the following types of indecomposable objects:  bands, bispecial strings,  finite and infinite special strings, finite and infinite 
strings.  We shall give more details in the overnext section. 

\begin{remark}
In fact we have shown (see \cite{Gelfand}) that the derived category of 
the Harish-Chandra modules over $SL_{2}(\mathbb R)$ is tame. 
\end{remark}

\section{Matrix problem for a general nodal algebra}

Let $A$ be a nodal algebra, which is supposed to be basic, $T$ its center
$\tilde{A} = \End_{A}({\rm rad}(A))$. Recall that we have 3 types of simple 
$A$-modules (see \cite{Pure}):
\begin{enumerate}
\item Such  simple left $A$-modules $U$ that $l_{A}(\tilde{A}\tens_{A} U) = 1$,
\item $l_{A}(\tilde{A}\tens_{A} U) = 2$, $l_{\tilde{A}}(\tilde{A}\tens_{A} U) = 2$.
\item $l_{A}(\tilde{A}\tens_{A} U) = 2$, $l_{\tilde{A}}(\tilde{A}\tens_{A} U) = 1$.
\end{enumerate}

It follows from the definition of a nodal algebra   that $A$ and $\tilde{A}$ have 
the common radical: $r = {\rm rad}(A) = {\rm rad}(\tilde{A})$. Hence we have an embedding
of semi-simple algebras $A/r \lar \tilde{A}/r$. Since $A$ is basic and
$\kk$ algebraically closed,
$A/r$ is  isomorphic to a product of several copies of $\kk$. 

The conditions 1-3 above imply that each simple component of 
$\tilde{A}/r$ is isomorphic either to $\kk$ or to $M_{2}(\kk)$ and the induced map
$A/r\lar \tilde{A}/r$ acts as follows:
\begin{enumerate}
\item A simple component of $A/r$ is mapped isomorphically onto a simple component of $\tilde{A}/r$.
\item A simple component of $A/r$ is embedded diagonally into a product of two simple components of $\tilde{A}/r$, both isomorphic to $\kk$.
\item A product of two simple components of $A/r$ is mapped isomorphically onto
the diagonal subalgebra of a simple component of $\tilde{A}/r$ isomorphic to
$M_{2}(\kk)$. 
\end{enumerate}

Let $I$ be an ideal in $\tilde{A}$ generated by the radical and idempotents of
the first type. Then $I$ is an ideal in $A$, too. Moreover, the factor-algebras
$A/I$ and $\tilde{A}/I$ are semi-simple in this case. So, the conditions of the
main theorem are fulfilled. 

Let $\tilde{A} = \prod\limits_{n=1}^{N} \tilde{A}_{n},$ where
all $\tilde{A}_{n}$ are hereditary orders, $C(\tilde{A}_{n})$ be the  basic algebra
corresponding to $\tilde{A}_{n}$. Since it is a hereditary order over
$\kk[[t]]$ (by Noether normalization there is a finite ring extension
$\kk[[t]] \lar T$), it is isomorphic to the completed path algebra of some
cycle of  length $d_{n}$ (it follows from the classification of hereditary orders over a complete discrete valuation ring, see \cite{Harada}, \cite{Brumer} or
\cite{DrozdKirich2}). Let us introduce some numbering
of the vertices of the cycles $C(\tilde{A}_{n})$. 
For the sake of  convenience we
number the vertices of $C(\tilde{A}_n)$ by elements  $[1], [2],\dots, [d_{n}]$
of ${\mathbb Z}/d_{n}{\mathbb Z}$.  So each simple
$\tilde{A}$-module $U$ correspond to a pair $(n,\nu)$, where $n\in 1,\dots, N, 
\nu\in {\mathbb Z}/d_{n}{\mathbb Z}.$ 
Namely, $n$ denotes the number of the component $\tilde{A}_n$ that acts non-trivially on $U$, $\nu$ is the number of the vertex from the cycle $C(\tilde{A}_n)$
corresponding to $U$.

In order to consider the category of triples ${\rm TC}_{A}$ we have to 
consider morphisms in the derived category $D^{-}(\tilde{A}-\mod)$. From 
what we have seen above it follows that  it is enough to consider morphisms
in $D^{-}(C(\tilde{A}_{n})-\mod), n\in 1,\dots, N$. 

Let $C$ be a cycle of length $m$. Then the category of finitely generated
 left
$C$-modules is hereditary. Hence any indecomposable object of
$D^{-}(C-\mod)$ is isomorphic to  $0\lar M \lar 0$, where $M$ is
an indecomposable $C$-module. 
Moreover, either $M$ is projective or it has a resolution 
 $P\stackrel{\varphi}\lar Q$, where $P$ and $Q$ are indecomposable projective $C$-modules, $\varphi$ a morphism, given by some path on the quiver $C$
\cite{uniserial}; denote $l(\varphi) = 
\mbox{\rm length}(\mbox{\rm coker}\varphi)$. 

The morphisms of $A$-modules (which are non-zero modulo the radical)
are of the following form (see \cite{uniserial}):

\begin{tabular}{p{1.7cm}ccc}
&
\xymatrix
{
P_{\nu+1} \ar[r]^{\varphi_{1}}  & P_{\nu} \\
P_{\nu+2} \ar[r]^{\varphi_{2}} \ar[u]^{c_{\nu+1}} & P_{\nu} \ar[u]\\
P_{\nu+3} \ar[r]^{\varphi_{3}} \ar[u]^{c_{\nu+2}} & P_{\nu} \ar[u]\\
\dots \ar[u]^{c_{\nu+3}} & \dots \ar[u]\\
}
&
\qquad
&
\xymatrix
{
P_{\nu} \ar[r]^{\psi_{1}} \ar[d] & P_{\nu-1}\ar[d]^{c_{\nu-1}}\\
P_{\nu} \ar[r]^{\psi_{2}} \ar[d] & P_{\nu-2} \ar[d]^{c_{\nu-2}}\\
P_{\nu} \ar[r]^{\psi_{3}} \ar[d] & P_{\nu-3} \ar[d]^{c_{\nu-3}}\\
\dots  & \dots, \\
}
\end{tabular}

\noindent
where $c_{\nu+i}: P_{\nu+i+1} \lar  P_{\nu+i}$ is the map given by
an arrow going from the vertex $\nu+i$ to $\nu+i+1$. 
There are also  morphisms in the derived category, which correspond to 
$\Ext$-groups:

\begin{tabular}{p{3.4cm}c}
&
\xymatrix
{
0 \ar[r] \ar[d]& P \ar[r] \ar[d]& Q  \ar[d]\\
0 \ar[r] \ar[d]& P \ar[r] \ar[d]\ & 0\ar[d]\\
R \ar[r] & P \ar[r] & 0\\
}
\end{tabular}

Let us now construct the partially ordered set, which describes the matrix problem corresponding to the category of triples ${\rm TC}_{A}$ for a given nodal algebra $A$.

Let $C(\tilde{A}_n)$ be a basic algebra (which is a cycle)
 corresponding to $\tilde{A}_n$.  Consider a complex
$$
(P_{\nu + l(\varphi)} \stackrel{\varphi}\lar P_{\nu})[f],
$$
where $\nu \in \ZZ/d_{n}\ZZ$, $\varphi$ a morphism of projective modules given 
by the path of the length  $l(\varphi)$, $f\in \ZZ$ the shift of the complexes.

Denote $\bJ(C(\tilde{A}_{n}))$ the set of simple $\tilde{A}_{n}$-modules, which correspond to direct summands of $\tilde{A}\tens_A U$, where $U$ is 
a simple $A$-module of second or third type. 

Let $\nu, \nu+l(\varphi) \in \bJ(C(\tilde{A}_{n}))$. Then we associate to this complex
 two symbols $\balpha(n,\nu,l(\varphi),f)$ and $\bbetha(n,\nu + l(\varphi),
l(\varphi), f+1)$.  In case only $\nu$ (respectively $\nu+l(\varphi)$ or neither of both)
belongs to $\bJ(C(\tilde{A}_{n}))$  we associate with it only
$\balpha(n,\nu,l(\varphi),f)$ (respectively $\bbetha(n,\nu + l(\varphi),
l(\varphi), f+1)$ or nothing). In the same way  a symbol
$\rho(n,\nu,f)$ corresponds to the object 
$$
(0\lar P_{\nu}\lar 0)[f]. 
$$

We are ready now to introduce our partially ordered set.

\begin{definition}
We introduce a Bondarenko's partially ordered set together with equivalence relation in several steps.
\begin{enumerate}
\item 
Let $1\le n \le N$, $\nu \in \ZZ/d_{n}\ZZ$, $\nu \in \bJ(C(\tilde{A}_n))$.
$$
\bE^{(f)}_{\nu}(n) = \{ \balpha(n,\nu,i,f), \bbetha(n,\nu,i,f)| i \ge 1\}
\cup\{\rho(n,\nu,f)\}, f\in \ZZ.
$$
\item 
$$
\bE_{\nu}(n) = \bigcup\limits_{f \in \ZZ}  \bE^{(f)}_{\nu}(n)
$$
$$
\bE(n) = \bigcup\limits_{\nu \in \bJ(\tilde{A}_n)}  \bE_{\nu}(n)
$$
$$
\bE = \bigcup\limits_{n=1}^{N}  \bE(n)
$$
\item 
We can introduce a partial order an an equivalence relation on $\bE$.
\begin{enumerate}
\item
First of all
$$
\balpha(n,\nu,i_{1},f) \ge \balpha(n,\nu,i_{2},f), \, \,
\bbetha(n,\nu,i_{1},f) \le \bbetha(n,\nu,i_{2},f)
$$
for $i_1 \ge i_{2}$. 
\item
Furthermore,
$$
\balpha(n,\nu,i,f) \ge \rho(n,\nu,f) \ge \bbetha(n,\nu,j,f)
$$
for all $i,j \ge 1$, $f\in \ZZ$. 
\item 
If $i\in {\mathbb N}$ and $\nu \in \bJ(C(\tilde{A}_n))$ are such that
$i+\nu \in \bJ(C(\tilde{A}_n))$, then 
$$
\balpha(n,\nu,i,f) \sim \bbetha(n,\nu + i,i,f+1), f\in \ZZ.
$$
\end{enumerate}
\item
Let $1\le n \le N$, $\nu \in \ZZ/d_{n}\ZZ, \nu \in \bJ(C(\tilde{A}_n))$. 
The set $\bF(n,\nu,f), f\in \ZZ$ consists either from one or two elements. 
\begin{enumerate}
\item
If $U$ is a simple module of a second type, 
$(n,\nu)$, $(m,\mu)$ corresponding  simple $\tilde A$-modules,
then the sets $\bF(n,\nu,f) = \{g(n,\nu,f)\}$ and 
$\bF(m,\mu,f) = \{g(m,\mu,f)\}$ are the  sets consisting from one element.
Moreover $g(n,\nu,f) \sim g(m,\mu,f)$. 
\item
In case $U$ is an $A$-module of the third type, $(n,\nu)$ corresponds
to \\ $B\tens_{A} U$, then $\bF(n,\nu,f) = \{g'(n,\nu,f), g''(n,\nu,f) \}$.
It is however convenient to assume that  
$\bF(n,\nu,f) = \{g(n,\nu,f) \}$ and $g(n,\nu,f)\sim g(n,\nu,f)$. 
\end{enumerate}
\end{enumerate}
\end{definition}

Let us point out that  we are interested in only in \emph{bounded from the right
representations} of the constructed bunch of semi-chains (since we want to describe the derived category of bounded from the right complexes). Moreover,
the non-degeneracy condition from the definition of the category of triples
implies certain non-degeneracy restrictions  on our matrices. However they 
concern only 
the discrete series of representations, for continuous series they are automatically satisfied (see \cite{mp}). 

To sum everything up we formulate the main result of this article:

\begin{theorem}
Let $A$ be a nodal algebra. The description of indecomposable objects
of $D^-(A-\mod)$ can be reduced to  the description of 
indecomposable representations of a  bunch of semi-chains, described in 
the previous definition. 
In particular, there
are 3 types of indecomposable objects in
$D^-(A-\mod)$:
\begin{enumerate}
\item Bands $\kB(w,m,\lambda)$.
\item Strings (which can be usual, special and  bispecial).
\item Complexes $P_i$ and $P_i \stackrel{\varphi}\lar P_j$, where $i$ and $j$
correspond to simple $A$-modules of the first type. 
\end{enumerate}
In particular,  a nodal algebra is derived-tame in ``pragmatic sense''.
\end{theorem}

\begin{remark}
For the ring $A= \kk[[x,y]]/(xy)$ it was was shown in \cite{curves}, 
how to describe complexes,
corresponding to objects of $A-\mod$ with respect to the canonical inclusion 
$$
A-\mod \lar D^{-}(A-\mod).
$$
\end{remark}

\section{Description of indecomposable complexes via gluing diagrams}

In this section we want to show, how the combinatoric of bands and strings can be  applied to write down explicit projective resolutions of indecomposable complexes.  We shall consider two ``typical examples'': the case of $A = \kk\langle\langle x,y\rangle\rangle/(x^2, y^2)$ and the case of the completed path algebra
of the Gelfand quiver. 

\subsection{The case of 
$D^{-}(\kk\langle\langle x,y\rangle\rangle/(x^2, y^2)$\mbox{\rm -- mod)}}

Let $A = \kk\langle\langle x,y\rangle\rangle/(x^2, y^2)$, $r$ its radical, 
$\tilde{A} = \mbox{\rm End}_{A}(r)$.   
As we have seen in previous sections, the description of indecomposable 
objects of the derived category  $D^{-}(A-\mod)$ can be  reduced to a matrix problem of type ``representations of bunches of chains''. There are two types of indecomposable complexes in this case: bands $\kB(w,m,\lambda)$ and 
strings $\kU(w)$.

Let us rewrite the corresponding partially ordered sets in this special case. 
We have a family of sets 
$\bF(k) = \{g(1,k), g(2,k)\}$, $k\in {\mathbb Z}$, $g(1,k)\sim g(2,k)$,  
which correspond to  columns  of the matrices $H_{k}(i|1)$ and $H_{k}(i|2)$;
and two families of sets 
$\bE_1(k) = \{ \alpha(1,i,k), \beta(1,i,k) | i\ge 1 \}\cup \{\rho(1,k) \}$ and 
$\bE_2(k) = \{ \alpha(2,i,k), \beta(2,i,k) | i\ge 1 \}\cup \{\rho(2,k) \}$, 
which label horizontal blocks of matrices 
$H_{k}(i|1)$ and $H_{k}(i|2)$, $k\in {\mathbb Z}$. 
 
The symbols $\{\rho(1,k) \}$ and $\{\rho(2,k) \}$ correspond to the $k$-th 
shift  of projective $\tilde{A}$-modules $P_1$ and $P_2$. 
The  element 
$
\beta(1+i,i,k+1)$ is conjugated to $\alpha(1,i,k)
$
and  
$
\beta(2+i,i,k+1)$ is conjugated to $\alpha(2,i,k),
$ 
where $1+i$ and $2+i$ have to be taken modulo 2. 

Let $w$ be some word containing a subword $\beta(1+i,i,k) \sim \alpha(1,i,k)$.
If $i$ is even, then it comes from the complex
$(P_1 \stackrel{\varphi}\lar  P_1)[k]$, where ${\rm coker}(\varphi)$ is an indecomposable $\tilde{A}$-module  of the 
length $i$. In what follows we shall say that $\varphi$ has length ${\rm coker}(\varphi)$.
 If $i$ is odd, then this subword  corresponds to 
$(P_2 \stackrel{\varphi}\lar P_1)[k]$.

As we shall see, an indecomposable complex from  the derived category 
$D^{-}(A-\mod)$ can be viewed as a gluing of complexes
$$
P_1, \, P_2, \, P_1 \stackrel{\varphi}\lar P_1, \, P_2 \stackrel{\psi}\lar P_2, \,
P_1 \stackrel{\phi}\lar P_2, \, P_2 \stackrel{\theta}\lar P_1.
$$
Suppose we have  a subword 
$$
 \beta(1+i,i,k+1) \sim \alpha(1,i,k) - g(1,k) \sim g(2,k) - \alpha(2,j,k) \sim
\beta(2+j,j,k+1). 
$$
It can be interpreted as a gluing  of complexes

\begin{tabular}{p{3.2cm}c}
&
\xymatrix
{ P_{i+1} \ar@{->}[r]^{\varphi_1} & P_1 \ar@{.}[d] \\
P_{j+2} \ar@{->}[r]^{\varphi_2} & P_2 \\
}
\end{tabular}

\noindent
shown  by the dotted line. Here the indices $i+1$  and $j+2$ must be 
taken modulo 2,  $\varphi_1$ and $\varphi_2$ have  the length $i$ and $j$ respectively.

\noindent
The subword 
$$
 \beta(1+i,i,k+1) \sim \alpha(1,i,k) - g(1,k) \sim g(2,k) - \beta(2+j,j,k) \sim
\alpha(2,j,k-1),
$$
 corresponds to the gluing of the type 

\begin{tabular}{p{2.8cm}c}
&
\xymatrix
{ P_{i+1} \ar@{->}[r]^{\varphi_1} & P_1 \ar@{.}[d]  & \\
& P_{j+2} \ar@{->}[r]^{\varphi_2} & P_2 \\
}
\end{tabular}

\noindent
and so on.

It is convenient to describe  gluing of the complexes by a gluing diagram.

\begin{example}
{\rm 
Consider the band data $\kB(w,m,\lambda)$, where 
$w = \alpha(2,2,0) \sim \beta(2,2,1) - g(2,1)\sim g(1,1) - \alpha(1,3,1) \sim \beta(2,3,2) - g(2,2)\sim g(1,2) - $ $\beta(1,2,2) \sim \alpha(1,2,1) - g(1,1)\sim g(2,1) -  \alpha(2,1,1) \sim \beta(1,1,0)-g(1,0)\sim g(2,0).$

Then it corresponds
to the following gluing diagram:

\begin{tabular}{p{3cm}c}
&
\xymatrix
{
P_2^m \ar@{->}[r]^{xyxI_m} \ar@{.}[ddd] & P_1^m \ar@{.}[d] &        \\
            & P_2^m \ar@{->}[r]^{xyI_m}  & P_2^m \ar@{.}[d]^{J_m(\lambda)} \\
                        & P_2^m \ar@{->}[r]^{xI_m} \ar@{.}[d]  & P_1^m  \\
P_1^m \ar@{->}[r]_{yxI_m}              & P_1^m                          &      \\
}
\end{tabular}

\noindent
This gluing diagram gives   the complex

\begin{tabular}{p{3cm}c}
&
\xymatrix
{
 &  A^m \ar@{->}[rd]^{xyI_m} & \\
A^m  \ar@{->}[ru]^{xyxI_m} \ar@{->}[rd]_{xyI_m} &  &  A^m \\
& A^m \ar@{->}[ru]_{xJ_m(\lambda)}&  \\
}
\end{tabular}

or, the same, 
$$
 A^m 
\xrightarrow{
\left(
\begin{array}{c}
xyxI_m \\
yxI_m \\
\end{array}
\right)}  A^{2m} 
\xrightarrow{
\left(
\begin{array}{cc}
xyI_m(\lambda) & xJ_m(\lambda) \\
\end{array}
\right)}  A^m.
$$
}
\end{example}

\begin{example}
{\rm 
Consider the string data $\kU(w)$, where 
$w = \dots - g(2,1)\sim g(1,1) - \beta(1,1,1) \sim \alpha(2,1,0) -g(2,0)\sim g(1,0) -\alpha(1,1,0) $ $ \sim  \beta(2,1,1)- g(2,1)\sim g(1,1) -
\alpha(1,1,1) \sim \beta(2,1,2) - \dots $ 

\noindent
The gluing diagram is

\begin{tabular}{p{2cm}c}
&
\xymatrix
{ P_2 \ar@{->}[r]^{x} \ar@{.}[d] & P_1 \ar@{.}[d] & P_2 \ar@{->}[r]^{x} \ar@{.}[d] & P_1 \ar@{.}[d] &  \\
 & P_2 \ar@{->}[r]^{x} & P_1 & P_2 \ar@{->}[r]^{x}  & P_1 \ar@{.}[d] \\
 & P_1 \ar@{->}[r]^{y} & P_2 & P_1 \ar@{->}[r]^{y}  & P_2  \\
P_1 \ar@{->}[r]^{x} \ar@{.}[u] & P_2 \ar@{.}[u] & P_1 \ar@{->}[r]^{x} \ar@{.}[u] & P_2 \ar@{.}[u] &  \\
}
\end{tabular}

\noindent
This string complex is  the minimal resolution  of 
 the simple module $\kk$
$$
\dots \lar A^2 
\xrightarrow{
\left(
\begin{array}{cc}
x & 0 \\
0 & y \\
\end{array}
\right)}  A^2 
\xrightarrow{
\left(
\begin{array}{cc}
x & 0 \\
0 & y \\
\end{array}
\right)}  A.
$$
}
\end{example}

There are also finite strings.

\begin{example}
{\rm 
\noindent
Consider the string $\kU(w)$ given by $$w = \rho(1,1)-g(1,1)\sim g(2,1) -\beta(2,3,1)\sim \alpha(1,3,0)- g(1,0)\sim g(2,0)-\rho(2,0).$$  

\begin{tabular}{p{3.8cm}c}
&
\xymatrix
{  P_1 \ar@{.}[d]& \\
   P_2 \ar@{->}[r]^{yxy} & P_1 \ar@{.}[d]\\
    & P_2 \\
}
\end{tabular}

\noindent
It corresponds to the complex
$$
 A \stackrel{yxy}\lar  A. 
$$
}
\end{example}

\subsection{The case of the Gelfand quiver}
Let $A$ be the completed path algebra of the Gelfand quiver, $r$ its radical and
$\tilde{A} = \mbox{\rm End}_{A}(r)$. 
Let $P_1, \, P_2,\,P_3,$ be indecomposable projective $A$-modules, $P,\,Q$ indecomposable projective $\tilde{A}$-modules and suppose that $\tilde{A}\tens_{A} P_1 = \tilde{A}\tens_{A} P_2 = P$ and $\tilde{A}\tens_{A} P_3 = Q$.

We have seen that the description of indecomposable objects of the derived category $D^{-}(A-$\mbox{\rm mod}) can be reduced to a matrix problem of the 
type ``representations of bunches of semi-chains''.

The combinatorics of indecomposable objects is similar to  the case of 
bunches of chains.  Continuous series of representations are still bands
$\kB(w,m,\lambda)$,
but the structure of discrete series is much more complicated. There are
 bispecial 
strings $\kU(w,m,\delta_1, \delta_2)$,  finite and infinite special strings 
$\kU(w,\delta)$  and  finite and infinite   usual strings $\kU(w)$. 
In this case there are also complexes (certain discrete series) which do not come from the matrix problem. 

Let us rewrite the partially ordered set in this case. We have a  family of 
sets 
${\mathbf F}(k) = \{ g(k)\}$ with equivalence relation  $g(k) \sim g(k)$
 The set ${\mathbf E}(k) = \{ \alpha(i,k), \beta(i,k)| i\in 
{\mathbb N} \} \cup \{\rho(k) \}$ is a chain with the total order
$$
\beta(j_2,k) \ge \beta(j_1,k) \ge \rho(k) \ge \alpha(i_1,k) \ge \alpha(i_2,k) 
$$
for all natural numbers $i_1\ge i_2$ and $j_1 \ge j_2$ and $k\in {\mathbb Z}$. 

If $i = 2l$ is even then we have conjugate points 
$\alpha(2l,k) $ and $ \beta(2l,k+1)$, and the subword  
$\alpha(2l,k) \sim \beta(2l,k+1)$ corresponds to the complex 
$(P \stackrel{\varphi}\lar P)[k],$ where  $\varphi$ is  the unique path from 
$b$ to itself  of the length $2l$. If 
$i = 2l+1$ is even, then  elements $\alpha(2l+1,k)$ and $\beta(2l+1,k)$
correspond to complexes $(P \stackrel{\varphi}\lar Q)[k]$ and 
$(Q \stackrel{\varphi}\lar P)[k-1]$ respectively, where $\varphi$ has the length $2l+1$.

As in the case of dihedral algebra, the combinatorics of bands and strings can be simplified. 
$$
\beta(2i,k+1) \sim \alpha(2i,k) - g(k) \sim g(k) - \alpha(2j,k) \sim
\beta(2j,k+1) 
$$
codes the gluing 

\begin{tabular}{p{4.5cm}c}
&
\xymatrix
{ P \ar[r]|{2i} & P \ar@{.}[d]\\
  P \ar[r]|{2j} & P \\
}
\end{tabular}

\noindent
etc. 
There is  an algorithm which associates to a band or string data the corresponding complex of projective modules. A complex of projective  $A$-modules is obtained as a gluing of the complexes
of $\tilde{A}$-modules $P \stackrel{\varphi}\lar P$, $Q \stackrel{\phi}\lar P $, $P \stackrel{\psi}\lar Q$ and $P$. In order to keep the notation simpler
we shall write instead of the map $\varphi$ only its length $l(\varphi)$ (which defines 
$\varphi$ uniquely).

\begin{example}
{\rm 
Consider the  band data: $(w,1,\lambda)$, where 
$
w = \alpha(2,0) \sim \beta(2,1) - g(1)\sim g(1) - \alpha(6,1) \sim \beta(6,2) - g(2)\sim g(2) - \beta(4,2)$ $\sim$
$\alpha(4,1)$ $ - g(1)\sim g(1) -  \beta(4,1) \sim \alpha(4,0)-g(0)\sim g(0).
$
It gives the following gluing diagram

\begin{tabular}{p{2.5cm}c}
&
\xymatrix
{ & P \ar@{.>}[d] \ar@{->}[r]|2 & P \\
 P  \ar@{->}[r]|4  & P & \\
 P \ar@{.>}[u] \ar@{->}[r]|6  & P & \\
  & P \ar@{.>}[u] \ar@{->}[r]|4 & P \ar@{.>}[uuu]|{\lambda}\\
}
\end{tabular}
}
\end{example}

Dotted lines here are directed: the direction of the arrow shows that there is 
a map of complexes which induces a non-zero map modulo the radical in the corresponding component of the complex. 

Now we introduce the rule ``of moving of an arrow'': 

\begin{enumerate}
\item 
Any time we have the situation 

\begin{tabular}{p{3.9cm}c}
&
\xymatrix
{ \bullet \ar@{.>}[d] \ar@{->}[r]|{a} & \bullet \\
\bullet &
}
\end{tabular}

we move the arrow (preserving its sign):

\begin{tabular}{p{3.9cm}c}
&
\xymatrix
{ \bullet \ar@{.>}[d] \ar@{->}[r]|{a} & \bullet \\
\bullet \ar@{->}[ur]|{a} &
}
\end{tabular}

\item 
Any time we have the situation 

\begin{tabular}{p{3.9cm}c}
&
\xymatrix
{   & \bullet  \ar@{.>}[d] \\
\bullet \ar@{->}[r]|{b} & \bullet
}
\end{tabular}

we can move the arrow taking the opposite sign:

\begin{tabular}{p{3.9cm}c}
&
\xymatrix
{   & \bullet  \ar@{.>}[d] \\
\bullet \ar@{->}[r]|{b} \ar@{->}[ur]|{-b}& \bullet
}
\end{tabular}

\item If we have the situation

\begin{tabular}{p{3.9cm}c}
&
\xymatrix
{ \bullet \ar@{.>}[d]|{\lambda} \ar@{->}[r]|{a} & \bullet \\
\bullet &
}
\end{tabular}

we move the arrow (preserving its sign):

\begin{tabular}{p{3.9cm}c}
&
\xymatrix
{ \bullet \ar@{.>}[d]|\lambda \ar@{->}[r]|{a} & \bullet \\
\bullet \ar@{->}[ur]|{\lambda a} &
}
\end{tabular}

\end{enumerate}

\noindent
Applying this rule to the band data above we get the following picture

\begin{tabular}{p{3cm}c}
&
\xymatrix
{ & P \ar@{.>}[d] \ar@{->}[r]|2 \ar@{->}[ddddr]|{\lambda 2}& P \\
 P  \ar@{->}[r]|6  \ar@{->}[ur]|6 \ar@{->}[ddr]|4 \ar@{->}[dddr]|{-4}& P \ar@{->}[ur]|{-2} \ar@{->}[dddr]|{-\lambda 2} & \\
 &    & \\
 P \ar@{.>}[uu] \ar@{->}[r]|4 \ar@{->}[dr]|{-4} & P \ar@{->}[dr]|{4}& \\
  & P \ar@{.>}[u] \ar@{->}[r]|4 & P \ar@{.>}[uuuu]|{\lambda}\\
}
\end{tabular}

Now we can insert instead of $P$ one of the symbols $P_1$ or $P_2$ following the rule that every dotted line has to connect symbols with different subscripts.

\begin{tabular}{p{3cm}c}
&
\xymatrix
{ & P_1  \ar@{->}[r]|2 \ar@{->}[ddddr]|{\lambda 2}& P_1 \\
 P_2  \ar@{->}[r]|6  \ar@{->}[ur]|6 \ar@{->}[ddr]|4 \ar@{->}[dddr]|{-4}& P_2 \ar@{->}[ur]|{-2} \ar@{->}[dddr]|{-\lambda 2} & \\
 &    & \\
 P_1  \ar@{->}[r]|4 \ar@{->}[dr]|{-4} & P_1 \ar@{->}[dr]|{4}& \\
  & P_2  \ar@{->}[r]|4 & P_2 \\
}
\end{tabular}

\noindent
It  corresponds to the  complex $\kP_\bullet$  
$$
P_2 \oplus P_1 \xrightarrow{
\left(
\begin{array}{rr}
\varphi_6 & 0 \\
\varphi_6 & 0 \\
\varphi_4 & \varphi_4  \\
-\varphi_4 & -\varphi_4 
\end{array}
\right)
} 
P_1 \oplus P_2 \oplus P_1 \oplus P_2 \xrightarrow{
\left(
\begin{array}{rrcc}
\varphi_2 & -\varphi_2 & 0 & 0 \\
\lambda \varphi_2 & -\lambda \varphi_2 &\varphi_4  & \varphi_4 
\end{array}
\right)
} 
P_1 \oplus P_2, 
$$
where 
$\varphi_{2k}$ always denotes  the map $\varphi_{2k}: P_i \lar P_j$ of the 
length $2k$.
Let us compute the triple $(\tilde{\kP}_\bullet, \kM_\bullet, i)$. Observe that 
after applying $\tilde{A}\tens_{A}$ to $P_i \stackrel{\varphi}\lar  P_j$ ($i,j \in \{1,2\})$  we get $P \stackrel{\varphi}\lar P$. It holds:
$$
\tilde{\kP}_\bullet := \tilde{A}\tens_{A} \kP_\bullet = \bigl(
P \oplus P \xrightarrow{
\left(
\begin{array}{rr}
\varphi_6 & 0 \\
\varphi_6 & 0 \\
\varphi_4 & \varphi_4  \\
-\varphi_4 & -\varphi_4 
\end{array}
\right)
} 
P \oplus P \oplus P\oplus P \xrightarrow{
\left(
\begin{array}{rrcc}
\varphi_2 & -\varphi_2 & 0 & 0 \\
\lambda \varphi_2 & -\lambda \varphi_2 &\varphi_4  & \varphi_4 
\end{array}
\right)
} 
P \oplus P \bigr). 
$$
Denote $M = M_2(\kk)$. The map $i: A/I\tens_{A} \kP_\bullet \lar 
\tilde{A}/I\tens \tilde\kP_\bullet$ is 

\begin{tabular}{p{0.4cm}c}
&
\xymatrix
{\kk_2 \oplus \kk_1   \ar@{->}[rr]^0 \ar[d]^{i_2} && \kk_1 \oplus \kk_2 \oplus 
\kk_1 \oplus \kk_2 \ar@{->}[rr]^0 \ar[d]^{i_1} && \kk_1 \oplus \kk_2 
\ar[d]^{i_0}\\
M\oplus M \ar@{->}[rr]^0 && M\oplus M\oplus M\oplus M \ar@{->}[rr]^0 && 
M\oplus M  \\
}
\end{tabular}

\noindent
where $i_0$, $i_1$ and $i_2$ are given by matrices

\begin{figure}[ht]
\hspace{1.3cm}
\includegraphics[height=3.5cm,width=9cm,angle=0]{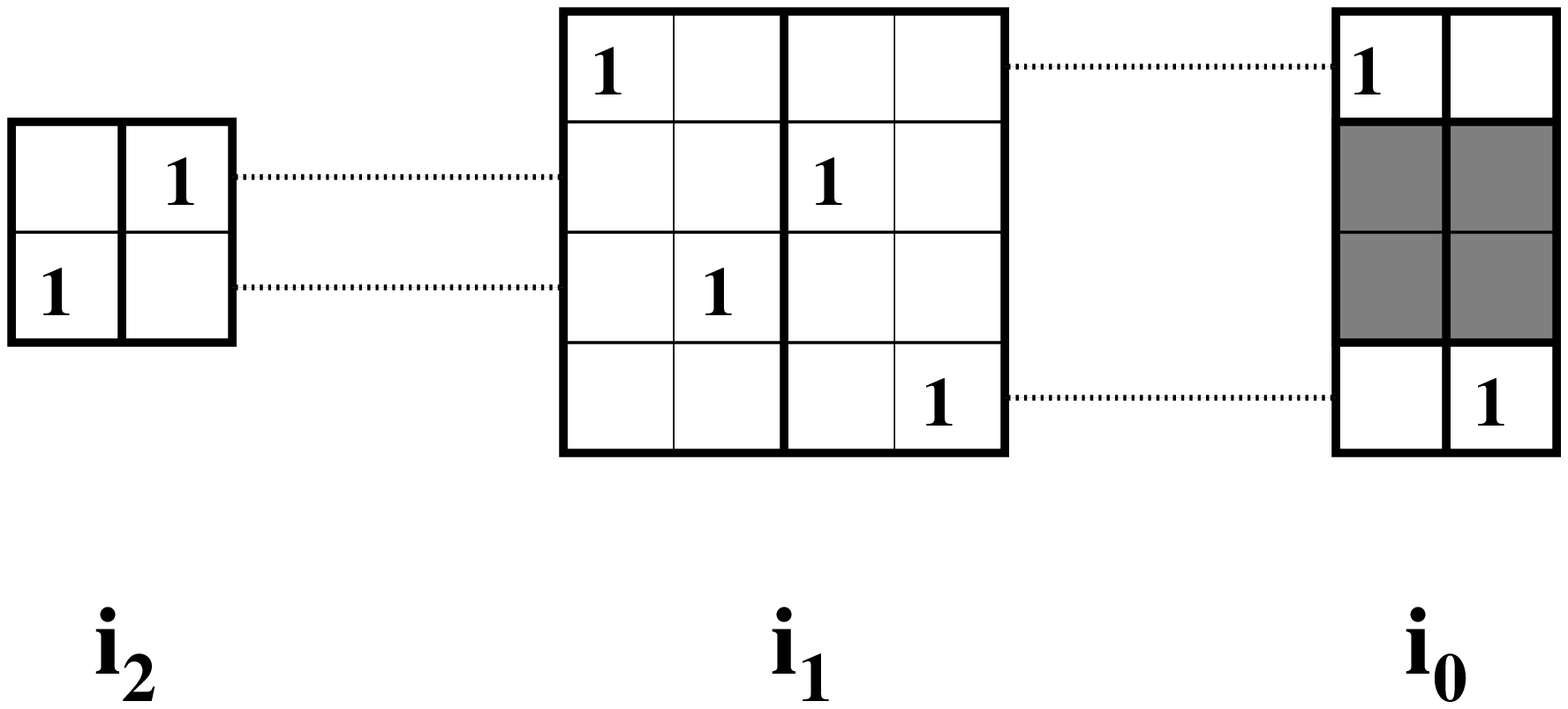}
\end{figure}

We have the following chain equivalence:

\noindent
\begin{tabular}{p{-3.5cm}c}
&
\xymatrix
{
P \oplus P  \ar[dddd]_{
\left(
\begin{array}{rc}
1 & 0 \\
1 & 1 \\
\end{array}
\right)
} \ar[rrr]^{
\left(
\begin{array}{rr}
\varphi_6 & 0 \\
\varphi_6 & 0 \\
\varphi_4 & \varphi_4  \\
-\varphi_4 & \varphi_4 
\end{array}
\right)
} 
& &&
P \oplus P \oplus P \oplus P \ar[dddd]^{
\left(
\begin{array}{rrrr}
1 & -1 & 0 & 0 \\
0 & 1 & 0 & 0 \\
0 & 0 & 1 & 0 \\
0 & 0 & 1 & 1 
\end{array}
\right)
}
\ar[rrr]^{
\left(
\begin{array}{rrcc}
\varphi_2 & -\varphi_2 & 0 & 0 \\
\lambda \varphi_2 & -\lambda \varphi_2 &\varphi_4  & \varphi_4 
\end{array}
\right)
} 
&&&
P \oplus P \ar[dddd]^{
\left(
\begin{array}{rc}
1 & 0 \\
-\lambda & 1 \\
\end{array}
\right)
}
\\
& & \\
& & \\
& & \\
P \oplus P \ar[rrr]^{
\left(
\begin{array}{cc}
0 & 0 \\
\varphi_6 & 0 \\
0 &  \varphi_4\\
0 & 0 
\end{array}
\right)
} 
& &&
P \oplus P \oplus P \oplus P \ar[rrr]^{
\left(
\begin{array}{rrcc}
\varphi_2 & 0 & 0 & 0 \\
0 & 0 & 0 & \varphi_4 
\end{array}
\right)
} 
&&&
P \oplus P 
}
\end{tabular}

\noindent
This map transforms the matrices $i_0$, $i_1$ and $i_2$ to the form

\begin{figure}[ht]
\hspace{1.3cm}
\includegraphics[height=3.5cm,width=9cm,angle=0]{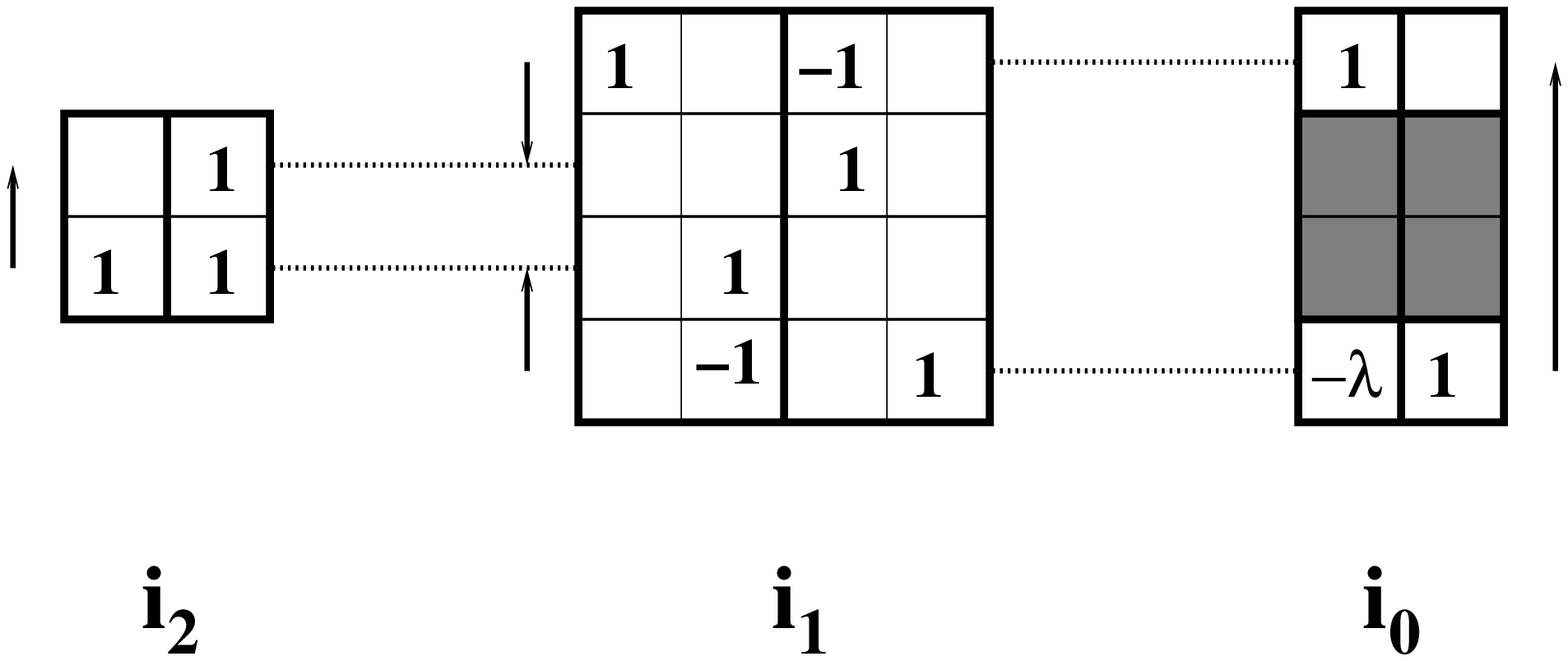}
\end{figure}

Doing the allowed transformations of rows and columns we transform them into the canonical form (see \cite{mp})

\clearpage

\begin{figure}[ht]
\hspace{1.3cm}
\includegraphics[height=3.5cm,width=9cm,angle=0]{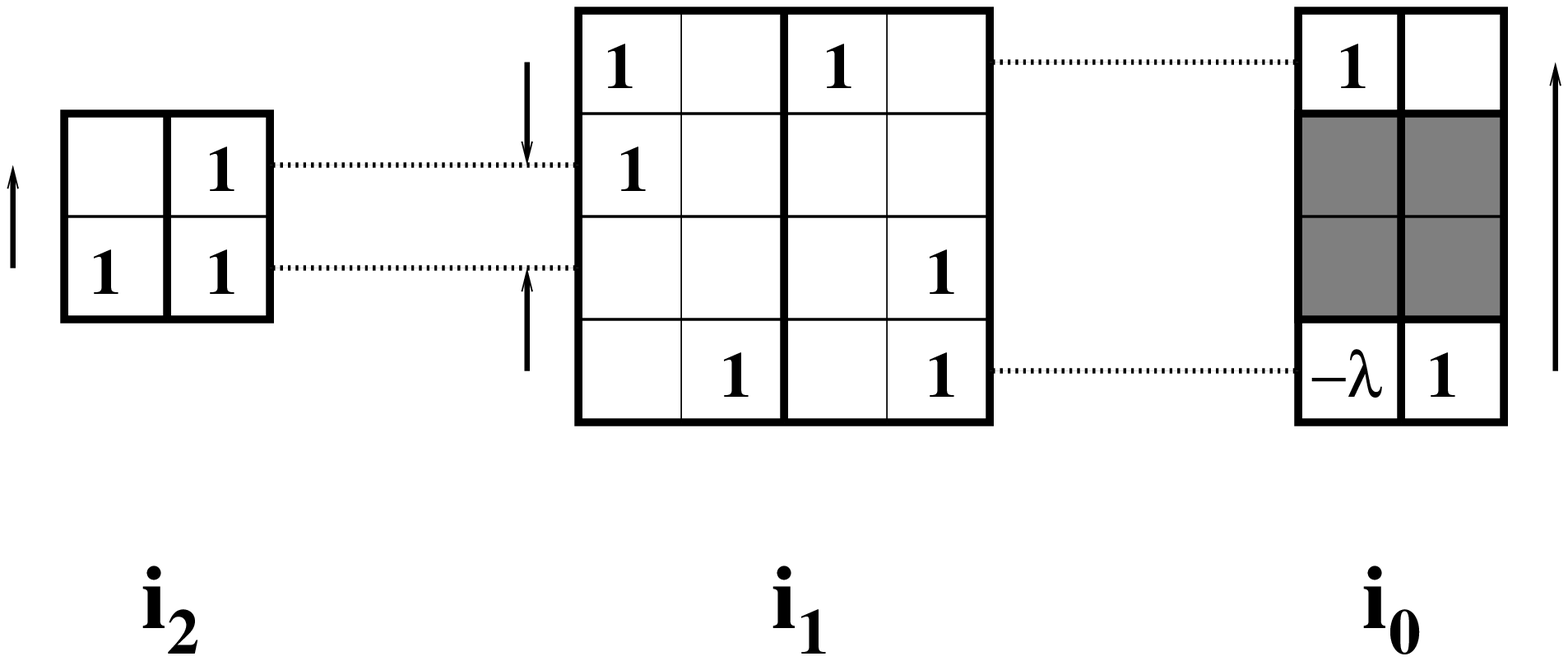}
\end{figure}

Suppose that a dotted line joins two points with equal weights. How to choose 
the direction of this line? We can do it by means of the following rule. 
Let us suppose that a gluing diagram has a subpart

\begin{tabular}{p{1.2cm}c}
&
\xymatrix
{   &                    &                                 & a  \ar@{.>}[d] \ar@{->}[r] & \bullet \\
    &                    &  \bullet \ar@{.>}[d] \ar@{->}[r]^3& \bullet               &             \\
 & x  \ar@{->}[r]^2& \bullet                        &                       &              \\
   \ar@{--}[rrrr]      &  &                               &                       &              \\
  & x  \ar@{->}[r]^2 \ar@{.>}[uu]  & \bullet           &                       &              \\
 &                       &  \bullet \ar@{.>}[u] \ar@{->}[r]^3  & \bullet  \ar@{.>}[d]             &             \\
 &                       &  \bullet  \ar@{->}[r]^3  & b               &             \\
}
\end{tabular}

We have to find first pair of points $(a,b)$ which are non-symmetric with respect to the axe of symmetry. 
In our case it   holds: $a < b$. The arrow looks in the direction of the  smaller point (see \cite{mp}).
In case when there are many dotted arrows joining points with equal weights, 
we have to consider  for each pair its own  axe of symmetry.

\clearpage  

\begin{example}
{\rm 
Consider the following gluing diagram (band):

\begin{tabular}{p{3cm}c}
&
\xymatrix
{          & \bullet \ar@{.>}[d] \ar@{->}[r]|2 & \bullet \ar@{.>}[ddddddd]|{\lambda}\\
\bullet  \ar@{->}[r]|2 & \bullet &                       \\
\bullet  \ar@{.>}[u] \ar@{->}[r]|2 & \bullet &                       \\
\bullet  \ar@{->}[r]|4 & \bullet \ar@{.>}[u] &                       \\ 
\bullet  \ar@{->}[r]|4 \ar@{.>}[u] & \bullet \ar@{.>}[d]&                       \\ 
\bullet  \ar@{->}[r]|2 \ar@{.>}[d]  & \bullet &                       \\ 
\bullet  \ar@{->}[r]|2              & \bullet &                       \\ 
         & \bullet \ar@{.>}[u] \ar@{->}[r]|6 & \bullet \\
}
\end{tabular}
}
\end{example}

If a word $w$ is symmetric, then we set directions of both dotted arrows intersecting the symmetry axe simultaneously clockwise or anticlockwise. 

\medskip

Let us now consider the case of discrete series. 
The first type of them are bispecial strings $\kU(w,m,\delta_1,\delta_2)$.
 They are given by some word
$w$, by a natural number $m$ and by two symbols 
$\delta_1, \delta_2 \in \{-,+\}$. 
Consider the following example: 

\begin{example}
{\rm 
Let $\delta_1 =+,\, \delta_2 =-$, $m=5$ and 
$w = (+) g(1) - \alpha(4,1)\sim \beta(4,2) - $ $ g(2)\sim g(2) $ $ - \beta(6,2)\sim $ $ \alpha(6,1) - g(1)\sim g(1) - \beta(2,1)\sim \alpha(2,0) - g(0) (-)$.  
Then we get the following gluing diagram 

\begin{tabular}{p{2.9cm}c}
&
\xymatrix
{  5P \ar@{.>}[d] \ar@{->}[r]^{4} & 5P (+)  & \\
   5P  \ar[r]^6 & 5P &  \\
  & 5P \ar@{.>}[u] \ar[r]^2 & 5P (-)\\
}
\end{tabular}

\noindent
It correspond to the  complex:
\clearpage

\begin{tabular}{p{0.9cm}c}
&
\xymatrix
{  5P_1  \ar@{->}[rrr]^{\varphi_4 I_5^{r+}} & & & 3P_1 \oplus 2P_2  & \\
   5P_2  \ar[rrr]^{\varphi_6 I_5}  \ar[drrr]_{- \varphi_6 I_5} \ar@{->}[urrr]^{\varphi_4 I_5^{r+}}& & & 5P_2 \ar[dr]^{\varphi_2 J_5^{r-}} & \\
& &   & 5P_1  \ar[r]^{\varphi_2 J_5^{r-} } & 2P_1 \oplus 3P_2 \\
}
\end{tabular}

\noindent
where $I_5^{r+}$ and $J_5^{r-}$ are the following matrices:
$$
I_5^{r+} = 
\left(
\begin{array}{rrrrr}
1 & 0 & 0 & 0 & 0 \\
0 & 0& 1& 0   & 0 \\
0 & 0& 0 & 0 & 1 \\
\hline
0 & 1 & 0 & 0 & 0 \\
0 & 0 & 0 & 1 & 0 \\
\end{array}
\right)
\qquad
J_{5}^{r-} = 
\left(
\begin{array}{rrrrr}
0 & 1 & 1 & 0 & 0 \\
0 & 0 & 0 & 1 & 1 \\
\hline
1 & 1 & 0 & 0 & 0 \\
0 & 0& 1& 1   & 0 \\
0 & 0& 0 & 0 & 1 \\
\end{array}
\right)
$$
The matrix $J_{5}^{r-}$ is  obtained by the following rule: we take the $5\times 5$ matrix
$$
J_5 = 
\left(
\begin{array}{rrrrr}
1 & 1 & 0 & 0 & 0 \\
0 & 1 & 1 & 0 & 0 \\
0 & 0& 1& 1   & 0 \\
0 & 0 & 0 & 1 & 1 \\
0 & 0& 0 & 0 & 1 \\
\end{array}
\right)
$$
and permute  rows putting rows  with even and  odd 
numbers in separate horizontal blocks. The superscript ``--'' means that the
block with even rows comes first.   The same rule applied to the
identity matrix $I_5$ gives    $I_5^{r+}$. 
} 
\end{example}

\noindent
The triple corresponding to this complex is isomorphic to  $(\tilde\kP_\bullet,
\kM_\bullet, i)$ where  $$\tilde\kP_\bullet = \bigl(
5P \oplus 5P \xrightarrow{
\left(
\begin{array}{cc}
I_5 & 0 \\
0   &I_5 \\
0 &  0 \\
\end{array}
\right)
}
5P \oplus 5P \oplus 5P \xrightarrow{
\left(
\begin{array}{ccc}
0 & 0 & I_5 \\
\end{array}
\right)
}
5P\bigr),
$$
$\kM_\bullet = \bigl(\kk_1^5\oplus \kk_2^5 \stackrel{0}\lar \kk_1^3\oplus \kk_2^2 \oplus \kk_1^5
\oplus \kk_2^5 \stackrel{0}\lar \kk_2^5 \oplus \kk_1^5\bigr) $
and  $i_0$, $i_1$, $i_2$ are given by matrices

\begin{figure}[ht]
\hspace{1.3cm}
\includegraphics[height=3.5cm,width=9cm,angle=0]{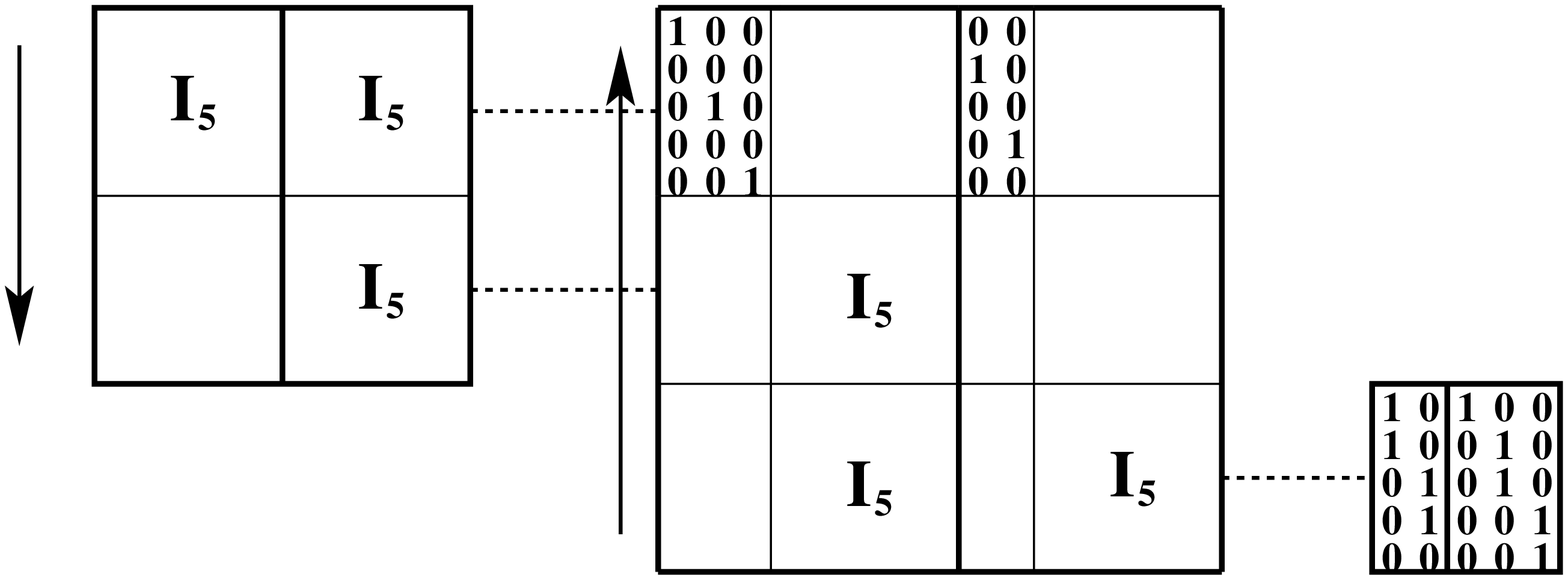}
\end{figure}

In the last  example both special ends 
were sinks. In the case when one of the special end is source we have to modify our rule a little bit.

\begin{example}
{\rm 
Consider the following bispecial string: $m = 4$, $\delta_1 =+, \, \delta_2 =+$,
$w = (+) g(2)-\beta(2,2)\sim $ $ \alpha(2,1) -g(1)\sim $ $ g(1) - \beta(2,1)\sim \alpha(2,0) - g(0) (+).$

\begin{tabular}{p{2.5cm}c}
&
\xymatrix
{ (+) P \ar[r]|2 &  P  & \\
                 &  P \ar@{.>}[u] \ar[r]|4 & P (+) \\
}
\end{tabular}

\noindent
It corresponds to the  complex

\begin{tabular}{p{0.7cm}c}
&
\xymatrix
{  2P_1 \oplus 2P_2 \ar[rr]^{\varphi_2 I_4^{c+}} \ar[drr]_{-\varphi_2 I_4^{c+}}& &  4P_1 \ar[drr]^{\varphi_4 J_4^{r+}} & &  \\
            & &  4P_2  \ar[rr]^{\varphi_4 J_4^{r+}} & & 2P_1 \oplus 2P_2 \\
}
\end{tabular}

\noindent
where 
$$
I_4^{c+} = 
\left(
\begin{array}{rr|rr}
1 & 0 & 0 & 0 \\
0 & 0 & 1 & 0 \\
0 & 1 & 0 & 0 \\
0 & 0 & 0 & 1\\
\end{array}
\right)
\quad \mbox{\rm and} \quad 
J_4^{r+} = 
\left(
\begin{array}{rrrr}
1 & 1 & 0 & 0  \\
0 & 0& 1& 1  \\
\hline
0 & 1 & 1 & 0  \\
 0& 0 & 0 & 1 \\
\end{array}
\right).
$$
The matrix $I_4^{c+}$ can be  computed by the following rule: we take the matrix
$$
I_4 = 
\left(
\begin{array}{rrrr}
1 & 0 & 0 & 0 \\
0 & 1 & 0 & 0 \\
0 & 0 & 1 & 0 \\
0 & 0 & 0 & 1\\
\end{array}
\right)
$$
and group  odd and even columns  into separate   blocks. 
}
\end{example}
The last remark concern the case when a  bispecial string is given by a 
``short'' diagram.

\begin{example}
{\rm 
Let $m = 5$, $\delta_1 =+, \delta_2 =- $  and $w$ is just 
$(+) P \stackrel{\varphi_2}\lar P (-)$. Then  the corresponding complex is 

\begin{tabular}{p{2.2cm}c}
&
\xymatrix
{ 3P_1 \oplus 2P_2  \ar[r]^{\varphi_2 C} &  2P_1 \oplus 3P_2 
}
\end{tabular}

\noindent
where  $C = I_5^{c+} \cdot (J_5^{r-})^{-1}$.
}
\end{example}

\noindent
There are also special and usual strings, which can be finite and infinite. 

\begin{example}
{\rm 
Let $w = \beta(1,1) - g(1)\sim g(1) - \alpha(1,1)$. Then
 the string $\kU(w)$ corresponds to the gluing diagram

\begin{tabular}{p{3.3cm}c}
&
\xymatrix
{                        & P \ar[r]|1 \ar@{.>}[d] & Q \\
Q \ar[r]|1    & p			  &   \\
}
\end{tabular}

\noindent
It defines  the complex 

\begin{tabular}{p{3.3cm}c}
&
\xymatrix
{                        & P_1 \ar[r]^{\varphi_1}  & P_3 \\
P_3 \ar[r]^{\varphi_1}  \ar[ru]^{-{\varphi_1}}   & P_2	\ar[ru]^{-{\varphi_1}}		  &   \\
}
\end{tabular}

\noindent
One can  recognize in this complex a projective resolution of the simple 
$A$-module $U_1$. 
}
\end{example}

\begin{example}
{\rm
Let $\delta = +$, $w = (+)g(0)- \alpha(2,0)\sim $ $ \beta(2,1) - g(1)\sim g(1) -
\alpha(4,1)\sim $ $ \beta(4,2) - $ $ g(2)\sim g(2) -
\alpha(2,2)$ $ \sim \beta(2,3) - g(3)\sim g(3) -\dots$.
The infinite special string $\kU(w,\delta)$ is given by the gluing diagram 

\begin{tabular}{p{1.6cm}c}
&
\xymatrix
{      \dots    &  P \ar[r]|2 \ar@{.}[d] & P         & P \ar[r]|2  \ar@{.}[d]& P (+) \\
\dots P \ar[r]|4 & P          & P\ar[r]|4 \ar@{.}[u]& P		   &  \\
}
\end{tabular}

\noindent
It corresponds to the complex 

\begin{tabular}{p{1.6cm}c}
&
\xymatrix
{      \dots    &  P_2 \ar[r]^{\varphi_2} \ar[ddr]^{-\varphi_2}& P_1 \ar[ddr]^{\varphi_4}        & P_2 \ar[r]^{\varphi_2}  & P_1 \\
 & &  & & \\
\dots P_2 \ar[r]^{\varphi_4} \ar[uur]^{-\varphi_4}& P_1 \ar[uur]^{\varphi_2}         & P_2\ar[r]^{\varphi_4} \ar[uur]^{-\varphi_4}	 & P_1 \ar[uur]^{\varphi_2}		   &  \\
}
\end{tabular}

\noindent
This complex belongs to $D^{-}(A$-\mbox{\rm mod}) and does not belong to 
$D^{b}(A$-\mbox{\rm mod})
}
\end{example}

\begin{example}
{\rm
Let $w = \rho(1,0)-g(0)\sim g(0)-\alpha(1,0)$. The the string complex $\kU(w)$ is given by the gluing diagram

\begin{tabular}{p{3.4cm}c}
&
\xymatrix
{ & P_1 \\
P_3 \ar[r]^{\varphi_1} \ar[ur]^{-\varphi_1}& P_2
}
\end{tabular}

\noindent
This complex is isomorphic to a projective resolution of a module which is finitely generated but not finite dimensional.
}
\end{example}

There are finally complexes which are not coming from the matrix problem. They are just complexes of the form $P_3 \stackrel{\varphi}\lar P_3$, which come from 
triples \\ $\bigl( Q \stackrel{\varphi}\lar Q, 0, 0 \bigr)$. 

The description of complexes for a general nodal algebra can be obtained by combining the combinatorics of complexes of the derived category of 
 the dihedral algebra and of the Gelfand quiver.

\section{Derived categories and  Harish-Chandra modules}
In \cite{Horoshkin} it was proven that there are only  two cases 
 of compact Lie groups, for which the category of Harish-Chandra modules 
is tame: 
$SL_{2}(\mathbb R)$ and $SO(1,n)$.  As a corollary of the theorem we obtain that
the derived category of Harish-Chandra modules is also tame in both of these cases.
We have already seen it for $SL_{2}(\mathbb R)$.  

Let $SO_{0}(1,n)$ be the connected component of $SO(1,n)$. 
\begin{enumerate}
\item Let $n=2l$. Then the category of Harish-Chandra modules over 
$SO_{0}(1,n)$ at a singular point
is equivalent to the category of finite-dimensional representations of the
completed path algebra of the following quiver (see \cite{Horoshkin}): 

\begin{figure}[ht]
\hspace{2.3cm}
\includegraphics[height=3.1cm,width=7.5cm,angle=0]{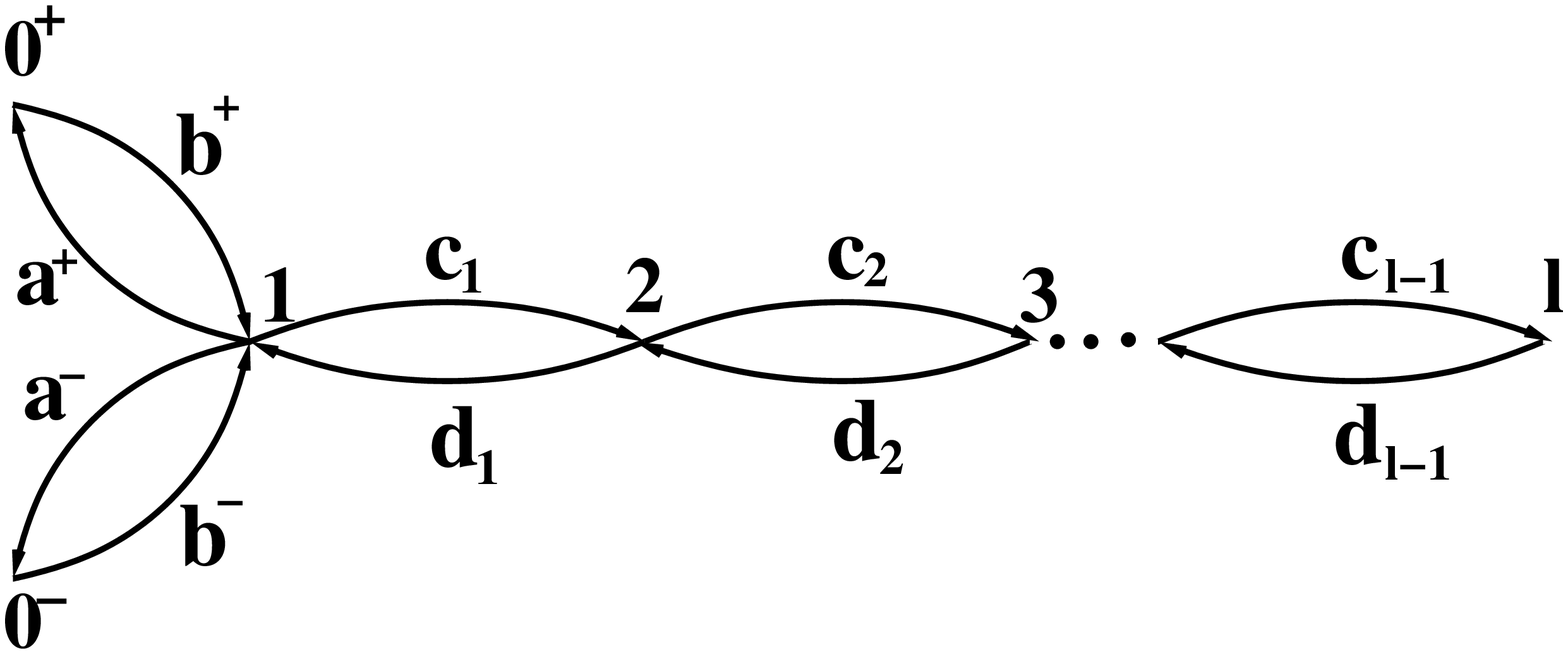}
\end{figure}

\noindent
where we have the relations:
$$
\gamma = b^{+}a^{+} = b^{-}a^{-}, 
$$
and $\gamma$ is  nilpotent,
$$
a^{\pm}d_{1} = 0, \hspace{0.2cm} c_{1}b^{\pm} = 0, 
$$
$$
c_{i+1} c_{i} = 0, \hspace{0.2cm} d_{i}d_{i+1} = 0, \hspace{0.2cm} i=1,
\dots,l-2.
$$
Moreover, all
$$
\vartheta_{i} = d_{i}c_{i}, \hspace{0.2cm} i=1,\dots,l-1.
$$
are nilpotent.

This algebra can be embedded into 

\begin{figure}[ht]
\hspace{1.6cm}
\includegraphics[height=2.8cm,width=9cm,angle=0]{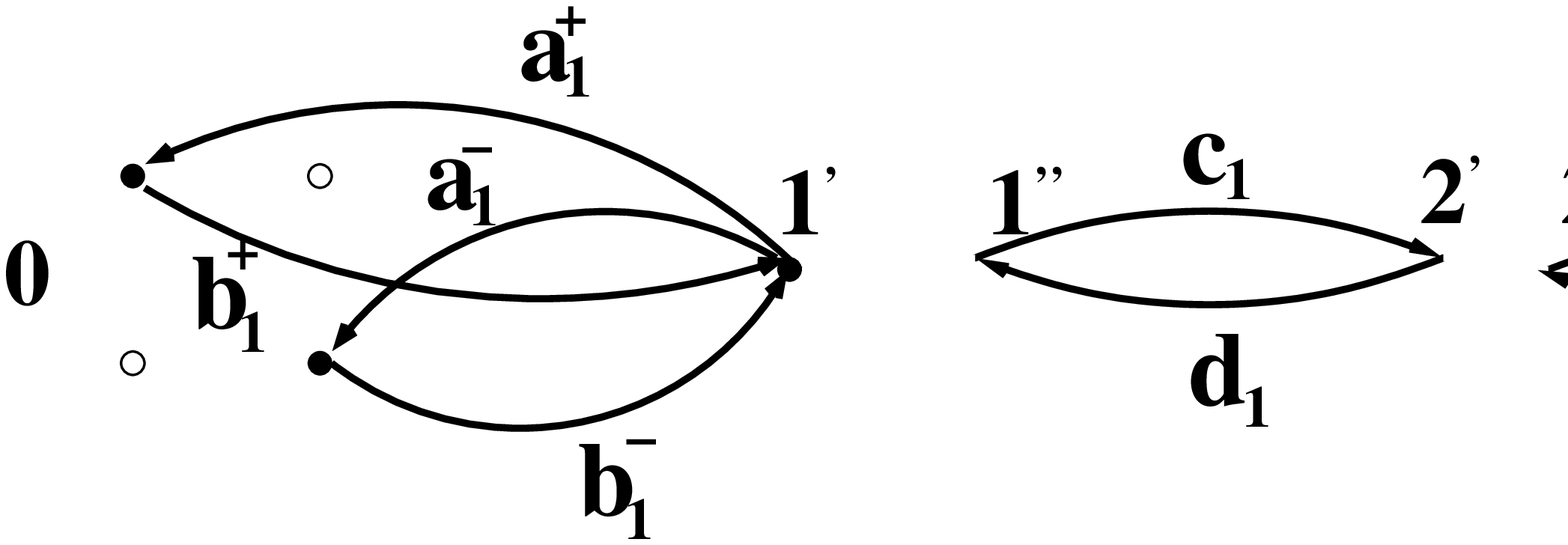}
\end{figure}

(this algebra is the endomorphism algebra of the radical of $A$. However,
it is not so important). 

The simple $A$-module, corresponding to the vertex $l+1$ is of the first type.
Those, which correspond to 
$1,2,\dots,l-1$ are of the second type, 
$0^{\pm}$ are of the third type.

\clearpage

\item If $n=2l+1$. Then the category of Harish-Chandra modules over 
$SO_{0}(1,n)$ is described by the completed path algebra of the 
following quiver:

\begin{figure}[ht]
\hspace{2.2cm}
\includegraphics[height=1.4cm,width=7cm,angle=0]{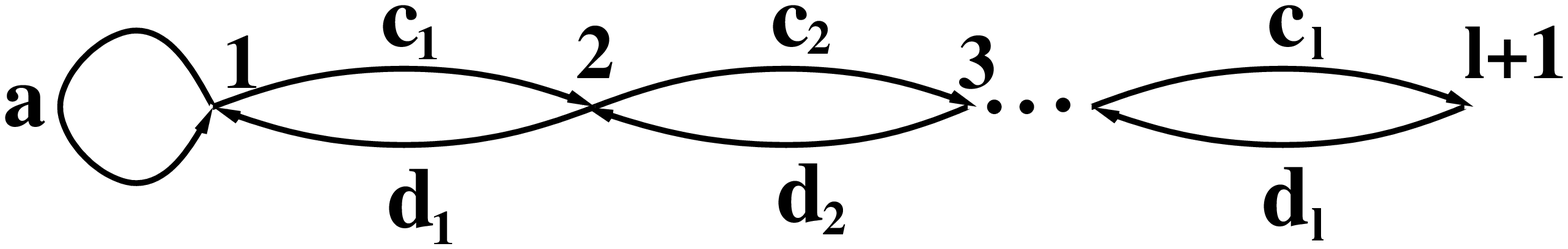}
\end{figure}

\noindent
with  relations:
$$
ad_{1} = 0, \hspace{0.2cm} c_{1}a = 0, 
$$
$$
d_{1}a = 0, \hspace{0.2cm} ac_{1} = 0, 
$$
$a$ is nilpotent,
$$
c_{i+1} c_{i} = 0, \hspace{0.2cm} d_{i}d_{i+1} = 0, \hspace{0.2cm} i=1,\dots,
l-1,
$$
and all
$$
\vartheta_{i} = d_{i}c_{i}, \hspace{0.2cm} i=1,\dots,l.
$$
are  nilpotent.

It can be embedded into

\begin{figure}[ht]
\hspace{1.2cm}
\includegraphics[height=1.6cm,width=8.5cm,angle=0]{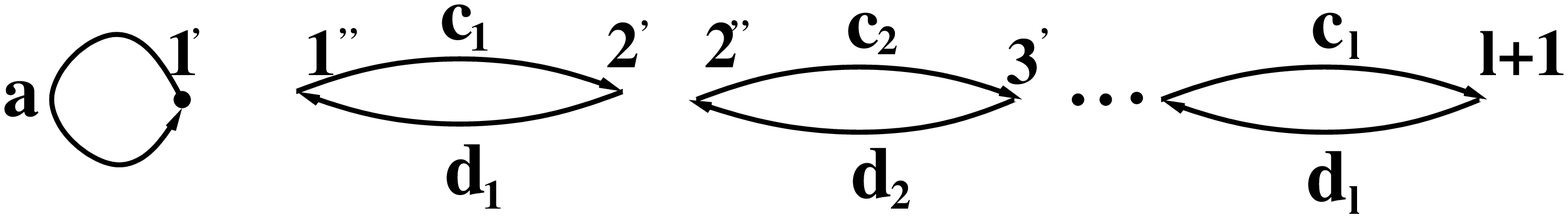}
\end{figure}

The simple module, corresponding to the vertex $l$ is of the first type,
all other simple modules   are of the second type. 
\end{enumerate}

Let us consider two more examples (see \cite{mp} for a description of
indecomposable modules over these algebras). 

\begin{example}
{\rm 
Consider the completed path algebra of the following quiver:

\begin{figure}[ht]
\hspace{1.7cm}
\includegraphics[height=3.5cm,width=9cm,angle=0]{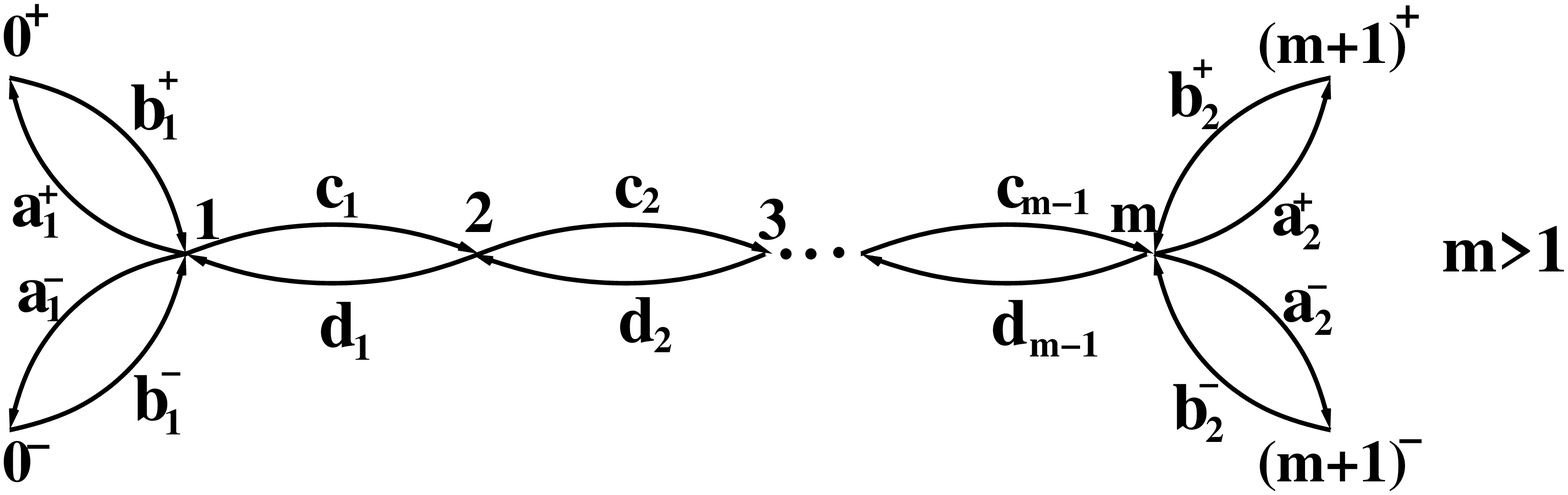}
\end{figure}

\noindent
where we have the relations
$$
b^{+}_{i}a^{+}_{i} = b^{-}_{i}a^{-}_{i}, \hspace{0.2cm} i=1,2,
$$
$$
a^{\sigma_{1}}_{i}b^{\sigma_{2}}_{i} = 0, \hspace{0.2cm}
 i=1,2, \quad \sigma_{1}, \sigma_{2} \in \{ -,+ \},
$$
$$
c_{i} d_{i} = 0, \quad  d_{i}c_{i} = 0, \hspace{0.2cm} 1\le i \le m-1
$$
and finally
$$
\gamma = b^{+}_{2}a^{+}_{2} c_{m-1} c_{m-2} \dots 
c_{1} b^{+}_{1}a^{+}_{1} d_{1}d_{2} \dots d_{m-1}
$$
is nilpotent.

As one easily observes, we can embed this algebra into  

\begin{figure}[ht]
\hspace{1.3cm}
\includegraphics[height=3.5cm,width=9cm,angle=0]{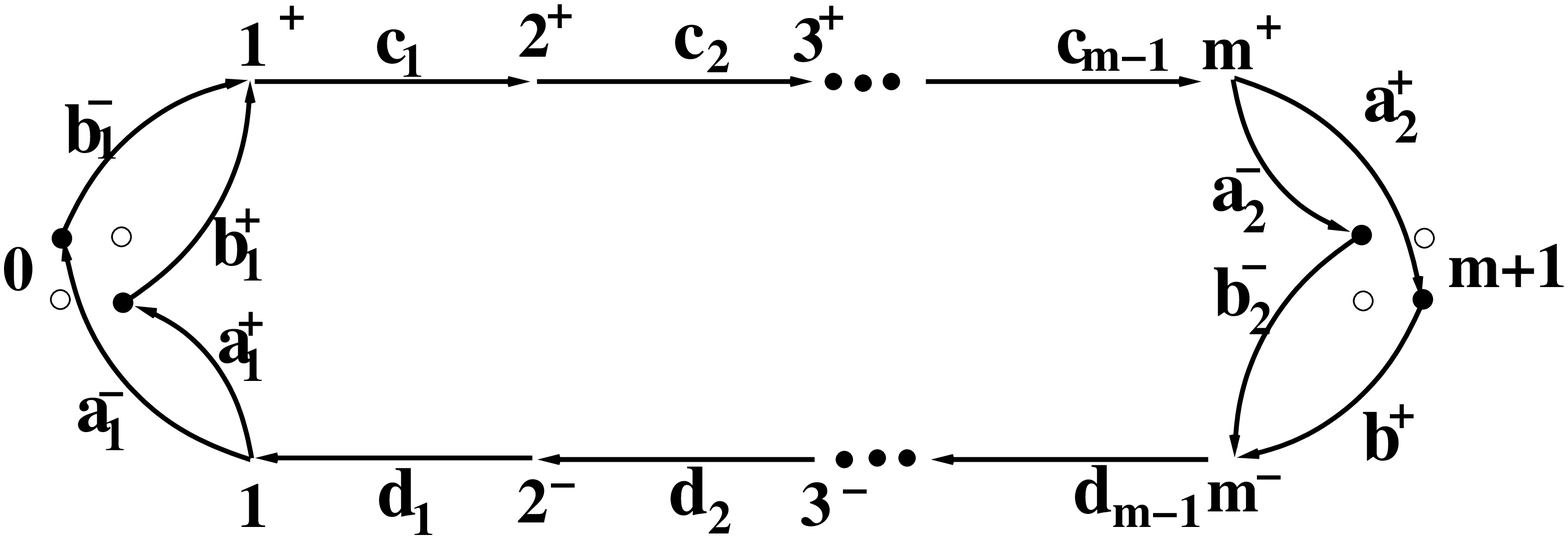}
\end{figure}

The simple $A$-modules, corresponding to vertices 
$1,2,\dots,m$ are of the second type, 
$0^{\pm}, (m+1)^{\pm}$ are of the third type.
}
 
\end{example}

\begin{example}
{\rm Consider the completed path algebra of quiver:

\begin{figure}[ht]
\hspace{1.5cm}
\includegraphics[height=3.5cm,width=9cm,angle=0]{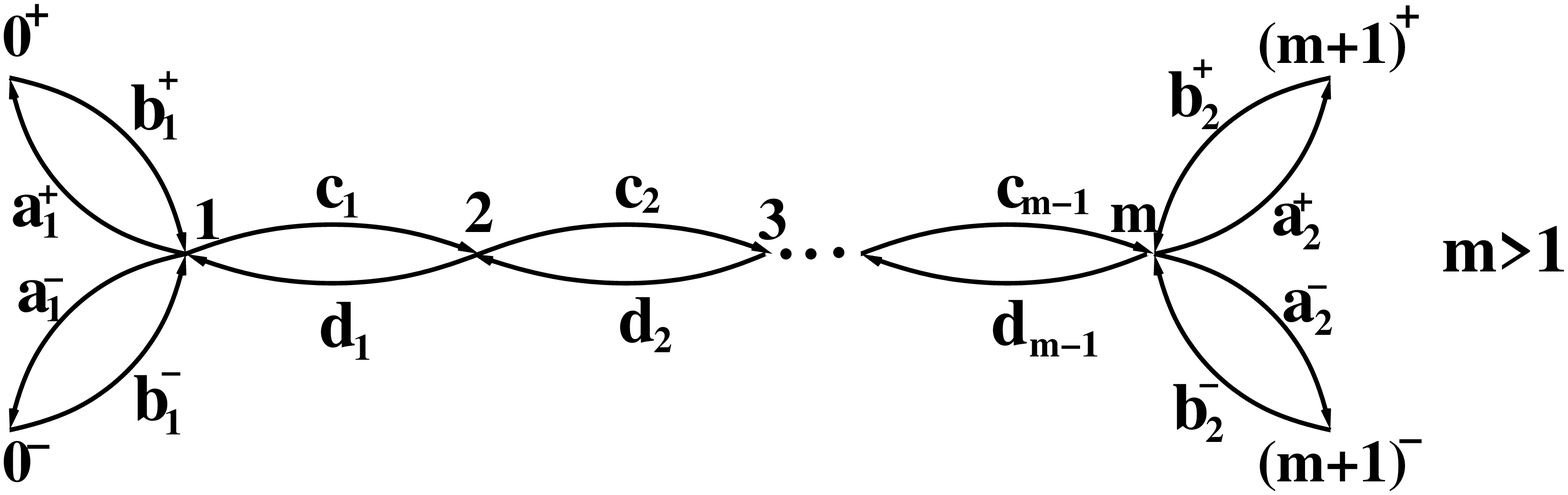}
\end{figure}

\noindent
where we have the relations:
$$
\gamma_{i} = b^{+}_{i}a^{+}_{i} = b^{-}_{i}a^{-}_{i}, \hspace{0.2cm} i=1,2,
$$
and $\gamma_{i}, i=1,2$ are nilpotent,
$$
a_{1}^{\pm}d_{1} = 0, c_{1}b_{1}^{\pm} = 0, 
$$
$$
a_{2}^{\pm}c_{m-1} = 0, d_{m-1}b_{2}^{\pm} = 0, 
$$
$$
c_{i+1} c_{i} = 0, d_{i}d_{i+1} = 0, \hspace{0.2cm} i=1,\dots,m-2.
$$
Moreover, all
$$
\vartheta_{i} = d_{i}c_{i}, \hspace{0.2cm} i=1,\dots,m-1.
$$
are nilpotent.

We can embed this algebra into:

\clearpage
\begin{figure}[ht]
\hspace{0cm}
\includegraphics[height=2.7cm,width=12cm,angle=0]{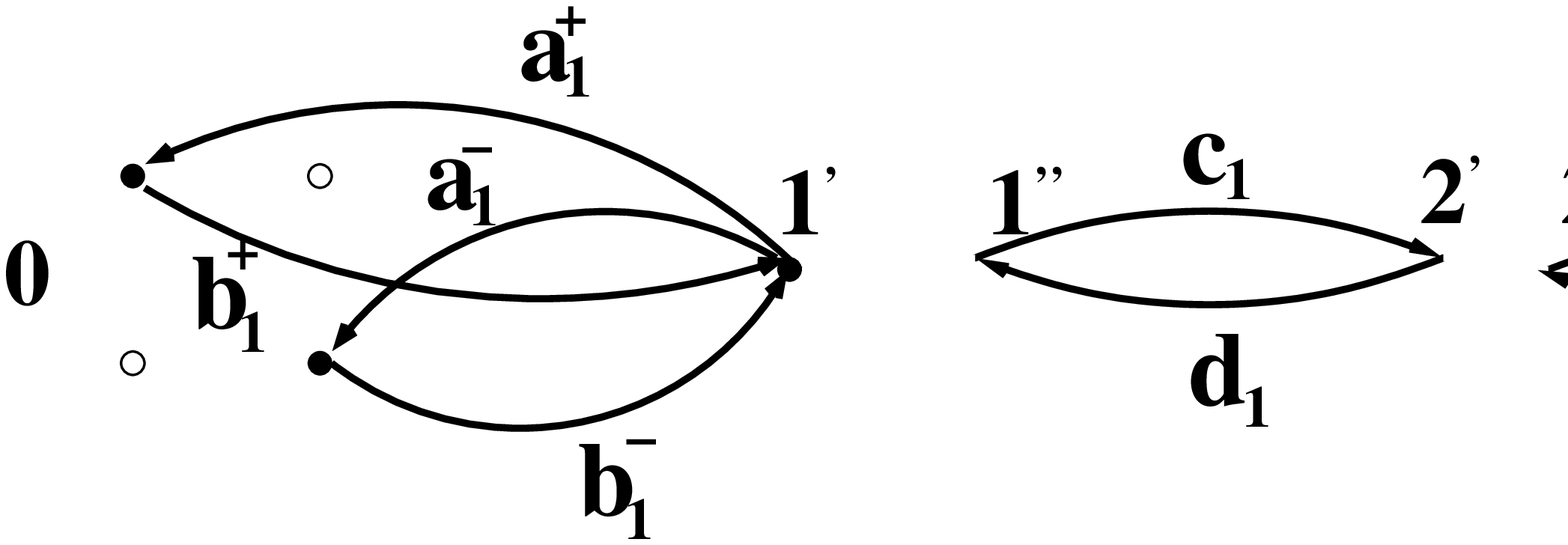}
\end{figure}

The simple $A$-modules, corresponding to vertices 
$1,2,\dots,m$ are of the second type, 
$0^{\pm}, (m+1)^{\pm}$ are of the third type.
}
\end{example}

\begin{remark}
It can be checked that all algebras from this section are nodal and  
every embedding is  embedding  into the endomorphism algebra of the radical. 
\end{remark}

\medskip
 \begin{center}
  \bf\Large Appendix
 \end{center} 

\appendix
\section{Krull--Schmidt theorem for homotopy categories}


\def\mZ{\mathbb Z}	\def\mN{\mathbb N}
\def\bL{\mathbf L}	\def\bS{\mathbf S}	\def\bR{\mathbf R}
\def\al{\alpha}	\def\be{\beta}	\def\fK{\mathbf k}
\def\im{\mathop\mathrm{Im}\nolimits}
\def\Mat{\mathop\mathrm{Mat}\nolimits}
\def\El{\mathop\mathrm{El}\nolimits}
\def\mtr#1{\begin{pmatrix}#1\end{pmatrix}}
\def\xx{\times}
\def\md#1#2#3{#1\equiv#2\pmod{#3}}

 \label{sec1a}

 Let $\,\dC\,$ be an additive category. We denote by $\,C(\dC)\,$ the category of  chain complexes with entries from $\,\dC\,$ and by $\,K(\dC)\,$ the factorcategory
 of $\,C(\dC)\,$ modulo homotopy. If $\,\dC\,$  is abelian, we denote by 
$\,D(\dC)\,$
 the \emph{derived category} of $\,\dC\,$, that is the category of quotients
 of $\,K(\dC)\,$ with respect to the set of morphisms inducing isomorphism of homologies.  We
 fix a commutative ring $\,\bS\,$ and consider $\bS$-\emph{categories}, namely, suppose that all sets
 $\,\dC(A,B)\,$ are modules over $\,\bS\,$ and the multiplication of morphisms is
 $\,\bS$-bilinear.

\begin{definition}\label{1a1}
 An additive category $\,\dC\,$ is called called:
\begin{itemize}
\item   \emph{local} if every object $\,A\in\dC\,$ decomposes into a finite direct sum
 of objects with local endomorphism rings;
\item  \emph{$\,\om$-local} if every object $\,A\in\dC\,$  decomposes into a finite of
 countable direct sum of objects with local endomorphism rings;
\item   \emph{fully additive}  if any idempotent morphism in $\,\dC\,$ splits, that is
 defines a decomposition into a direct sum;
 \item
  \emph{\lof} (over $\bS$) if it is fully additive and all morphism spaces $\dC(A,B)$ are
 finitely generated $\bS$-modules. Especially if $\bS$ is a field, a \lof category is
 called \emph{locally finite dimensional}.
\end{itemize}
\end{definition}
 Evidently, every \lof category is local; moreover, an endomorphism algebra $\dC(A,A)$
 in a \lof category is a \emph{finite $\bS$-algebra}, i.e. such that the underlying $\bS$-module
 is finitely generated. It is known that any local (or $\,\om$-local) category is fully additive; moreover, a
 decomposition into a direct sum of objects with local endomorphism rings is always unique;
 in other words, any local (or $\,\om$-local) category is a \emph{Krull--Schmidt} one,
 cf.~\cite[Theorem 3.6]{bass}.

 For a local category $\,\dC\,$ denote  by $\,\rad\dC\,$ its \emph{radical}, that is
 the set of all morphisms $\,f:A\to B\,$, where $\,A,B\in\ob\dC\,$, such that no component
 of the matrix presentation of $\,f\,$ with respect to some (hence any) decomposition
 of $\,A\,$ and $\,B\,$ into a direct sum of indecomposable objects
 is invertible. Note that if $\,f\notin\rad\dC\,$, there is a morphism $\,g:B\to A\,$
 such that $\,fgf=f\,$ and $\,gfg=g\,$. Hence both $\,gf\,$ and $\,fg\,$ are nonzero
 idempotents, which define decompositions $\,A\iso A_1\+A_2\,$ and $\,B\iso B_1\+B_2\,$
 such that the matrix presentation of $\,f\,$ with respect
 to these decompositions is diagonal: $\,\mtr{f_1&0\\0&f_2}\,$, and $\,f_1\,$ is invertible.
 Obviously, if $\,\dC\,$ is locally finite dimensional, then $\,\rad\dC(A,B)\,$ coincide with the set of all
 morphisms $\,f:A\to B\,$ such that $\,gf\,$ (or $\,fg\,$) is
 nilpotent for any morphism $\,g:B\to A\,$.
   
\begin{proposition}\label{1a2}
 Suppose that $\bS$ is a complete local noetherian ring with maximal ideal $\gM$.
 If $\,\dC\,$ is a  \lof category over $\bS$, the categories $\,C(\dC)\,$ and $\,K(\dC)\,$
  are $\,\om$-local (in particular, Krull--Schmidt).
 Moreover, a morphism $\,f\dt:A\dt\to B\dt\,$ from one of these
 categories belongs to the radical \iff all components $\,f_ng_n\,$ (or $\,g_nf_n\,$)
 are nilpotent modulo $\gM$ for any morphism $\,g\dt:B\dt\to A\dt\,$.
\end{proposition}
{\sc Proof.}\ 
 Denote by $\fK=\bS/\gM$ the residue field of $\bS$. We use the following simple statement from linear algebra.

 \begin{lemma}\label{12}
 Let $\,\La\,$ be a finite dimensional $\,\fK$-algebra and $\,a\,$ be and element from
 $\,\La\,$. There is a polynomial $\,f(x)\in\fK[x]\,$ such that $\,f(a)\,$ is an idempotent
 and $\,f(e)=e\,$ for any idempotent $\,e\,$ from any $\fK$-algebra. Moreover, 
 $\,f(a)\,$ is nilpotent (or invertible) \iff so is $\,a\,$.
\end{lemma}

\noindent
{\sc Proof.}\ 
Suppose that a polynomial $f(x)$ satisfies the condition $f(0)=0$, $f(1)=1$.
Then $f(e)=e$ for any idempotent $e$ from any finite dimensional algebra.
  
 We can embed $\,\La\,$ in an endomorphism algebra of some finite dimensional vector space
 $\,V$, so we suppose that $\,\La=\End V\,$. Decompose $\,V=V_0\+V_1\,$ so that the
 restriction $\,a|_{V_0}\,$ is nilpotent and $\,a|_{V_1}\,$
 is invertible. Replacing $\,a\,$ by $\,a^k\,$ for some $\,k\,$, one can suppose that
 $\,a|_{V_0}=0$. Indeed, if we have found a polynomial 
$f(x)$ such that $f(a^k)$ is idempotent, then  $f^{k}(a) = f(a^k)$ hence
$f^k(x)$ is the polynomial for $a^k$ we are looking for. In particular, if 
$a^k = 0 $, then we can take $f(x)=x^k$.  
 Set $\,b=a|_{V_1}\,$
 Since $\,b\,$ is invertible, there is a polynomial $\,g(x)\,$ such that $\,g(b)=1\,$
 and $\,g(0)=0\,$. If $\,1\,$ is an eigenvalue
 of $\,b\,$, then $\,g(1)=1\,$, whence $\,g(e)=e\,$ 
for every idempotent $\,e\,$.
 If $\,1\,$ is not an eigenvalue of $\,b\,$,
 then $\,(x h(x),x^2-x)=x\,$, where $\,h(x)\,$ is the minimal polynomial of $\,b\,$,
 hence there is a polynomial $\,f(x)\,$ such that
 $\,f(x)\equiv g(x)\pmod{xh(x)}\,$ and $\,f(x)\equiv x\pmod{x^2-x}\,$.
 Therefore,
 $\,f(b)=1\,$ and $\,f(e)=e\,$ for every
 idempotent $\,e\,$, which accomplishes the proof of the lemma.

\hfill$\Box$

 Recall also a known result, which can be easily deduced, for instance, from \cite[Section III.8]{jac}.
   \begin{lemma}
  There are polynomials $G_n(x)\in\mZ[x]$ with $G_n(0)=0$ and such that for every ring $\La$, any ideal $I\subseteq\La$
 and any element $a\in\La$ such that $a^2\equiv a\mod I$, $G_n(a)^2\equiv G_n(a)\mod I^{n+1}$ and $a\equiv G_n(a)\mod I$.
 \end{lemma} 
 (For instance, $G_1(x)=3x^2-2x^3$.)
 
  \begin{corollary}\label{ca1}
 Let $\La$ be a finite algebra over a local noetherian ring $\bS$ with maximal ideal $\gM$
 and $a\in\La$. For every $n\in\mN$ there is a polynomial $g(x)\in\bS[x]$ such that
 \begin{itemize}
\item   $g(a)^2\equiv g(a)\mod \gM^{n+1}$;
 \item
  $g(e)\equiv e\mod\gM^n$ for every element $e$ of an arbitrary finite $\bS$-algebra 
 such that $e^2\equiv e\mod\gM^n$;
 \item
  $g(a)\equiv1\mod\gM$ \iff $a$ is invertible;
 \item
  $g(a)\equiv0\mod\gM$ \iff $a$ is nilpotent modulo $\gM$.
\end{itemize}
 \end{corollary} 
 {\sc Proof.}
  Find a polynomial $f(x)$ over $\fK=\bS/\gM$ such that $f(\bar a)$ is idempotent in $\La/\gM\La$, where
 $\bar a=a+\gM\La\in\La/\gM\La$ and
 $f(\bar e)=\bar e$ for any idempotent $\bar e$ of any $\fK$-algebra; especially $f(0)=0$. Lift $f(x)$ to a polynomial
 $F(x)\in\bS[x]$ such that $F(0)=0$. Then $F(a)$ is idempotent modulo $\gM$ and if $e^2\equiv e\mod\gM^n$,
 then $F(e)\equiv e\mod\gM$ (by the construction  of $F(x)$) and $eF(e)\equiv F(e)\mod\gM^n$ (it is true for any polynomial $F(x)$ satisfying $F(0) = 0$). Set $g(x)=G_n(F(x))$. Then $g(a)$ is idempotent
 modulo $\gM^{n+1}$, just as $g(e)$ for every $e$ that is idempotent modulo $\gM$. If, moreover, $e^2\equiv e\mod\gM^n$,
 then $g(e)\equiv e\mod\gM$ and $eg(e)\equiv g(e)\mod\gM^n$. Let $g(e)=e+r$; then $r = g(e)-e$ and $er = re = r$. Therefore it holds 
 $(e+r)^2\equiv e+2r+r^2\equiv e+r
 \mod\gM^n$, wherefrom $r\equiv -r^2 \mod\gM^n$. But then 
$r\equiv -r^2 \equiv  -r^4 \equiv\dots \mod\gM^n$, so $r\equiv 0\mod\gM^n$. 
 \hfill$\Box$

 Let now $\,a\dt\,$ be an endomorphism of a complex $\,A\dt\,$ from $\,C(\dC)\,$.
 Consider the sets $\,I_n\subset\mZ\,$ defined as follows: $\,I_0=\set{0}\,$,
 $\,I_{2k}=\setsuch{l\in\mZ}{-k\le l\le k}\,$ and $\,I_{2k-1}=\setsuch{l\in\mZ}{-k<l\le k}\,$.
 Obviously,  $\,\bup_nI_n=\mZ\,$, $\,I_n\subset I_{n+1}\,$ and $\,I_{n+1}\=I_n\,$
 consists of a unique element $\,l_n$. Using corollary \ref{ca1}, we can construct a sequence of
 endomorphisms $a^{(n)}\dt$ such that
 \begin{itemize}
\item   $(a^{(n)}_i)^2\equiv a^{(n)}_i\mod\gM^n$;
 \item
  $a^{(n+1)}_i\equiv a^{(n)}_i\mod\gM^n$;
 \item
  $a^{(n)}_i$ is invertible or nilpotent modulo $\gM$ \iff so is $a_i$.
\end{itemize}
 Then one easily sees that setting $u_i=\lim_{n\to\8}a^{(n)}_i$, we get an idempotent endomorphism $u\dt$ of $A\dt$,
 such that $u_i\equiv0 \mod\gM \,\, (u_i \equiv 1\mod\gM)$ \iff $a_i$ is nilpotent modulo $\gM$ (respectively $a_i$ is invertible).

 Especially, if either one of $\,a_l\,$ is neither nilpotent nor invertible modulo $\gM$ or one of $\,a_l\,$ is nilpotent modulo $\gM$
 while another one is invertible, then  $\,u\dt\,$ is neither zero nor identity. Hence the complex $\,A\dt\,$ decomposes.
 Thus $\,A\dt\,$ is indecomposable if and only if, for any endomorphism $\,a\dt\,$ of $\,A\dt\,$,
 either $\,a\dt\,$ is invertible or all components $\,a_n\,$ are nilpotent modulo $\gM$.
 Since all algebras $\,\End A_n/\gM\End A_n\,$ are finite dimensional, neither
 product $\,\al\be,$ where $\,\al,\be\in\End A_n\,$ and one of them is nilpotent
 modulo $\gM$, can be invertible. Therefore, the set of endomorphisms
 $\,a\dt\,$ of an indecomposable complex $\,A\dt\,$ such that all components $\,a_n\,$
 are nilpotent modulo $\gM$ form an ideal $R$ of $\,\End A\dt\,$ and $\End A\dt/R$ is a skew field.
 Hence $R=\rad(\End A\dt)$ and $\,\End A\dt\,$ is local.

Now we want to show that any complex from $\,C(\dC)\,$ has an indecomposable 
 direct summand.  Consider an arbitrary complex $\,A\dt\,$ and suppose that 
$A_0 \ne 0$. For any idempotent endomorphism  $e\dt$ of $A\dt$
at least one of the complexes $e(A\dt)$ or $(1-e)(A\dt)$  has a non-zero component  at the zero place. On the set of all endomorphisms of $\,A\dt\,$ we 
can introduce  a partial ordering by writing 
 $e\dt\ge e'\dt$ \iff $e'\dt=e\dt e'\dt e\dt$ and both $e_0$ and $e'_0$ 
are non-zero. Let  $e\dt\ge e'\dt\ge e''\dt\ge\dots$ be a chain 
 of idempotent endomorphisms of $A\dt$. As all endomorphism algebras $\,\End A_l\,$ are finitely
 generated $\bS$-modules, the sequences 
$e_l,e'_l,e''_l,\dots \in \End A_l $ stabilize for all $l$, so this chain has a
 lower bound (formed by the limit values of components). By Zorn's lemma, there is a minimal non-zero idempotent of $A\dt$,
 which defines an indecomposable direct summand.

 Again, since all $\End A_l$ are finitely generated,
 for every $n$ there is a decomposition
 $\,A\dt=B^{(n)}\dt\+\bigoplus_{i=1}^{r_n}{B_{in}}\dt\,$ where all $\,{B_{in}}\dt\,$ are indecomposable
 and $\,B^{(n)}_l=0\,$ for $\,l\in I_n\,$. Moreover, one may suppose that $\,r_n\le r_m\,$
 for $\,m>n\,$ and $\,{B_{in}}\dt={B_{im}}\dt,$ for $\,i\le r_n\,$. Evidently, it implies
 that $\,A\dt=\bigoplus_{i=1}^r{B_i}\dt\,$ where 
 $\,r=\sup_nr_n\,$ and $\,{B_i}\dt={B_{in}}\dt\,$ for $\,i\le r_n\,$, which accomplishes
 the proof of the Proposition \ref{1a2} for $\,C(\dC)\,$. 

 Note now that the endomorphism ring of each complex $\,{B_i}\dt\,$ in the category
 $\,K(\dC)\,$ is a factor-ring of
 its endomorphism ring in $\,C(\dC)\,$. Hence it is either local or zero; in the latter
 case the image of $\,{B_i}\dt\,$ in $\,K(\dC)\,$ is a zero object. Therefore, the claim
 is also valid for $\,K(\dC)\,$. 
\hfill$\Box$

 \begin{corollary}\label{1a3}
  Let $\bS$ be local, complete and noetherian, and $\bA$ be an $\bS$-algebra finitely
 generated as $\bS$-module. Then the derived category $D^-(\mdl\bA)$, where $\mdl\bA$ denotes
 the category of finitely generated $\bA$-modules, is $\om$-local, in particular, Krull--Schmidt category.
 \end{corollary} 
 {\sc Proof.}\ 
   Indeed, $D^-(\mdl\bA)$ coincide with the category $K^-(\pro\bA)$, where $\pro\bA$ denotes the category
 of finitely generated projective $\bA$-modules.
 \hfill$\Box$ 

\begin{remark}\em
 The conditions of Proposition \ref{1a2} are essential indeed, and Krull--Schmidt
 theorem can fail even for the category of bounded
 complexes $\,C^b(\dC)\,$ over a local category $\dC$ as the following example shows.

 Let $\,\bR\,$ be the localization of the polynomial ring $\,\fK[x,y]\,$ at the maximal
 ideal $\,(x,y)\,$,
 $\,\dC=\pro\bR\,$ be the category of free $\,\bR$-modules of finite rank. Obviously, it is local.
 The category $\,K^-(\dC)\,$
 is equivalent to the category $\,D^-(\mdl\bR)\,$ and contains the category $\,\mdl\bR\,$ as a full subcategory.
 Denote by $\,\bS\,$ the factor-ring $\,\bR/(x^2y-y^3+x^4)\,$. It is a
 local domain, but its completion $\,\hat\bS\,$ is not a domain: its normalization 
 decomposes as $\,\bS_1\xx\bS_2\xx\bS_3\,$,
 where each $\,\bS_i\iso\fK[[x]]\,$. In particular, $\,\hat\bS\,$ has three torsion-free
 modules $\,\bL_i\,$ such that each $\,\bL_i\,$ has a composition series with the factors
 $\,\bS_j,\bS_k\,$, where $\,\set{i,j,k}=\set{1,2,3}\,$ (it is the projection of $\hat\bS$.
 onto $\bS_j\xx\bS_k$). 
 It implies that $\,\bS\,$ has torsion-free indecomposable modules
 $\,M_1,M_2,N_1,N_2,N_3\,$ with the following completions:
\begin{align*}
 \hat M_1&=\bS_1\+\bS_2\+\bS_3;\\
 \hat M_2&=L_1\+L_2\+L_3;\\
 \hat N_i&=\bS_i\+L_i\quad (i=1,2,3).
\end{align*}
  (cf.~\cite{rog}). Then $\,M_1\+M_2\iso N_1\+N_2\+N_3\,$,
 hence the category $\,\bS\mdl\,$, thus also $\,\bR\mdl\,$
 and $\,D^b(\bR$--mod)\, are not local.
\end{remark}

\section{Bunches of chains}
\label{seca}

We summarize the results of \cite{mp,mp1}, changing both the definition and the presentation of the answer to
 equivalent ones, which seem more convenient for our purpose. Moreover we add a description of some
 infinite dimensional representations that occur in dealing with derived categories, together with a sketch of
 a proof. As V.M.Bondarenko has informed us, he has submitted a paper containing more details on infinite case \cite{Bond}.  Note that Bondarenko calls ``bunch of semichained sets'' what
 we call ``bunch of chains''. The reason can be seen if one compares our definitions.

 \begin{definition}\label{a11}
  A \emph{bunch of chains} $\dX$ consists of:
 \begin{itemize}
\item  An index set $I$, which we suppose finite or countable.
 \item
  For each $i\in I$, two chains (linearly ordered sets) $\dE_i$ and $\dF_i$.\\
 We set $\dE=\bup_{i\in I}\dE_i,\ \dF=\bup_{i\in I}\dF_i$ and $|\dX|=\dE\cup\dF$.
 \item
  A symmetric relation $\sim$ (not an equivalence!) on $|\dX|$ such that for every $x$ there is at most one $y$
 with $x\sim y$ (maybe $x=y$).
\end{itemize}
 We define an equivalence relation $\approx$ on $|\dX|$ such that $x\approx y$ means either $x=y$ or $x\sim y$,
 and set $\ti\dX=|\dX|/\approx$.
 We write $x-y$ if there is an index $i\in I$ such that $x\in\dE_i,\,y\in\dF_i$ or vice versa. For each
 $x\in|\dX|$ such that $x\sim x$ we introduce two new elements $x',\,x''$ and
 set $\dX^*=(|\dX|\=\setsuch{x}{x\sim x}) \cup\setsuch{x',x''}{x\sim x}$. We subdivide $\dX^*$ into
 $\dE^*=\bup_i\dE_i^*$ and $\dF^*=\bup_i\dF^*_i$, which are the images of $\dE_i$ and $\dF_i$;
 for instance $x'$ and $x''$ are in $\dE^*_i$ if $x\in\dE_i$. We consider the ordering
 $<$ on $|\dX|$, which is just the union of orderings on all $\dE_i$ and $\dF_i$, and extend it, as well as
 the relation $-$, onto $\dX^*$ so that each ``new'' element $x'$ or $x''$ inherits all relations that the element
 $x$ has. For instance, $x'<y$ with $y\in|\dX|$ means that $x<y$; $x''-z'$ means that $x-z$, etc.
 Note that the elements $x',\,x''$ are always non-comparable. On the other hand, we extend the equivalence $\approx$ to
 $\dX^*$ trivially (each new element $x'$ or $x''$ is unique in its $\approx$-class), and set $\ti\dX^*=\dX^*/\approx$.
 \end{definition} 

 A bunch of chains $\dX$ gives rise to a bimodule problem. Namely, we fix a field $\fK$
 and define a $\fK$-category $\bA=\bA(\dX)$ and an $\bA$-bimodule $\bU=\bU(\dX)$ as follows:
 \begin{itemize}
 \item   
 $\ob\bA=\ti\dX^*$.
 \item
    If $a,b $ are two equivalence classes, a basis of the
  morphism space $\bA(a,b) $ consists of elements $p_{yx} $ with
  $x\in a,\, y\in b,\, x<y $ and, if $a=b $, the identity
  morphism $1_x $.
 \item
  The multiplication is given by the rule: $p_{zy}p_{yx}=p_{zx} $
  if $z<y<x $, while all other possible products are zeros.
 \item
  A basis of $\bU(a,b) $ consists of elements $u_{yx} $ with
  $y\in b\cap\dE^*, \,\,x\in a\cap\dF^*,\,x-y $.
 \item
  The action of $\bA $ on $\bU $ is given by the rule:
  $p_{zy}u_{yx}=u_{zx} $ if $y<z $; $u_{yx}p_{xt}=u_{yt} $
  if $x<t $, while all other possible products are zeros.
 \end{itemize}

 The category of \emph{representations of the bunch} $\dX $ over the field $\fK$ is then
 defined as the category $\El(\bU) $ of the elements of this
 bimodule. In other words, a representation is a set $M$ of block matrices 
$$ 
   M_i=\mtr{\hdotsfor{3}\\\dots&M_{xy}&\dots\\\hdotsfor{3}},
  \  i\in\bI,\ x\in\dE^*_i,\ y\in\dF^*_i,\ M_{xy}\in\Mat(n_x\xx n_y,\fK)
$$ 
 such that $x\approx y$ implies $n_x=n_y$. Two representations are isomorphic
 \iff they can be obtained from one another by a sequence of the following 
 \emph{elementary transformations}:
 \begin{itemize}
 \item 
 elementary transformations of rows (columns) in each horizontal (vertical) stripe;
 it means that they are performed simultaneously in all matrices $M_{xy}$
 with fixed $x$ (respectively $y$); moreover, if $x\approx z$, the transformations of the $x$-stripe
 must be the same as those of $z$-stripe (certainly, if one of them is horizontal and the other is vertical,
 ``the same'' means ``contragredient'');
 \item
 if $x<y$, then scalar multiples of rows (columns) of the $x$-stripe can be added to rows
 (columns) of the $y$-stripe.
 \end{itemize}
 One easily sees that this definition coincides with that of \cite{mp,mp1}.

 The description of indecomposable representations from \cite{mp,mp1} rests upon a combinatorics,
 which we expound in terms of \emph{strings and bands} alike to their use in the representation theory.

\begin{definition}\label{a12}
 Let $\dX=\set{I,\dE_i,\dF_i,\sim} $ be a bunch of chains.
\begin{enumerate}
 \item
  An $\dX $-{\em word} is a sequence $w=x_1r_1x_2r_2x_3\dots r_{m-1}x_m$,
 where $x_k\in|\dX| $ and $r_k\in\set{\sim,-} $, such that for all possible values of $k$
  \begin{enumerate}
   \item
    $x_kr_kx_{k+1} $ in $|\dX| $.
   \item
    $r_k\ne r_{k+1}$.
  \end{enumerate}
 We call $m$ the \emph{length} of the word $w$.
  Possibly $m=1 $, i.e. $w=x $ for some $x\in|\dX| $. The elements $x_1$ and $x_m$ are called
 the \emph{ends} of the word $w$.
 \item
  We call an $\dX $-word {\em full\/} if, whenever $x_1 $ is not a
  unique element in its $\approx$-class, then $r_1=\sim $, and
  whenever $x_m $ is not a unique element in its $\approx$-class, then $r_{m-1}=\sim $.
 \item
  We denote by $w^*$ the \emph{inverse word} $x_mr_{m-1}x_{m-1}\dots r_1x_1$ and
 call an $\dX$-word \emph{symmetric} if $w= w^*$. We call $w$ \emph{quasisymmetric}
 if it can be presented in the form $v\sim v^*\sim v\sim v^*\sim\dots\sim v$ for a shorter word $v$.
 \item
  We call the end $x_1$ ($x_m$) of the word $w$ \emph{special} if $x_1\sim x_1$ and $r_1={\!}-$
 (respectively $x_m\sim x_m$ and $r_{m-1}={\!}-$). We call the word $w$
   \begin{itemize}
  \item
   \emph{usual} if neither of its ends is special;
    \item
   \emph{special} if one of its ends is special;
  \item
   \emph{bispecial} if both its ends are special.
   \end{itemize}
\emph{ Note that a special word is never symmetric, while a bispecial word is always full; a quasisymmetric word
 is always bispecial.}
   \item
  If $r_1=r_{m-1}=\sim$ and $x_m-x_1$ in $\dX$,
 we call the word $w $ an {\em $\dX $-cycle}. 
  Note that in this case $m $ is always even. For a cycle $w$ we set $r_m=-$ and
 $x_{qm+k}=x_k,\,r_{qm+k}=r_k$ for all integers $q,k$.
 \item
  We call an $\dX $-cycle $w=x_1r_1x_2r_2x_3\dots r_{m-1}x_m $ {\em
  non-periodic} if  the sequence $x_1r_1x_2r_2\dots x_mr_m $
  cannot be written as a multiple $vv\dots v $ of a shorter sequence $v $.
 \item
  A \emph{shift} of a cycle $w$ is defined as the cycle $w^{[k]}=x_{k+1}r_{k+1}x_{k+2}\dots r_{k-1}x_k$ 
 for some even integer $0\le k< m$. We call a non-periodic cycle $w$ \emph{symmetric} if $w^*= w^{[k]}$ for some $k$. 
 (Note that $w^{[k]}=w^{[l]}$ with $k\ne l$ is impossible if $w$ is non-periodic.)
 \item
  For a cycle $w$ and an even integer $0\le k<m$ we define $\nu(k,w)$ as the number of even integers $0\le i \le k$ such that
 both $x_{i-1}$ and $x_i$ belong either to $\dE$ or to $\dF$.
\end{enumerate}
\end{definition}

  \begin{definition}\label{a13}
  \begin{enumerate}
 \item 
  A \emph{usual string} datum is a non-symmetric full usual word.
 \item
  A \emph{special string} datum is a pair $(w,\de)$, where $w$ is a special full word and $\de\in\set{+,-}$.
 \item
  A \emph{bispecial string} datum is a quadruple $(w,m,\de_1,\de_2)$, where $w$ is a bispecial word, which is
 neither symmetric nor quasisymmetric, $m\in\mN$ and $\de_i\in\set{+,-}$.
 \item
  A \emph{band} datum is a pair $(w,f)$, where $w$ is a non-periodic cycle and $f\in\fK[t]$ is a \emph{primary
 polynomial} over the field $\fK$, i.e. a degree of an irreducible polynomial with leading coefficient $1$,
 such that $f(0)\ne0$ and if $w$ is symmetric also $f(1)\ne0$. 
 If the field $\fK$ is algebraically closed and $f=(t-\la)^d$, we write $(w,d,\la)$ instead of $(w,f)$. 
 \item
 The following string data are called \emph{equivalent}:
  \begin{enumerate}
  \item 
   usual string data $w$ and $w^*$;
 \item
  special string data $(w,\de)$ and $(w^*,\de)$;
 \item
   bispecial string data $(w,m,\de_1,\de_2)$ and $(w^*,m,\de_2,\de_1)$.
  \end{enumerate}
 \item
  Two band data are called \emph{equivalent} if they can be obtained from one another by a sequence of the following transformations:
 \begin{enumerate}
 \item
  replace $(w,f)$ by $(w^{[k]},f)$ if $\nu(k,w)$ is even;
 \item
  replace $(w,f)$ by $(w^{[k]},\al^{-1}t^df(1/t))$, where $d=\deg f$ and 
$\al=f(0)$ if $\nu(k,w)$ is odd;
\item   replace $(w,f)$ by $(w^*, f)$. 
\end{enumerate}
 \emph{Note that if $f(t)=(t-\la)^d$, then $\al^{-1}t^df(1/t)=(t-\la^{-1})^d$.}
 \end{enumerate}
 \end{definition} 
 
 Then the main result of the papers \cite{mp,mp1} (see also \cite{CB}) can be reformulated as follows.

  \begin{theorem}\label{a14}
  There is one-to-one correspondence  between isomorphism classes of indecomposable representations of a bunch of chains and equivalence
 classes of string and band data.
 \end{theorem} 

 We call indecomposable representations corresponding to usual string data (special string data, bispecial string data,
 band data) \emph{usual strings} (respectively, \emph{special strings, bispecial strings, bands}). 

 \section{Infinite chains}

  For our purpose we have to  consider some infinite representations of a bunch of chains. We suppose
 now that $I=\mN$ and for every index $i\in I$ the set 
 $$
 \setsuch{j\in I}{(\exists x\in\dE_i\cup\dF_i)\, (\exists y\in\dE_j\cup\dF_j)\ x\sim y}
 $$
  is finite. Namely, we define the category $\El^\8(\bU)$ just in the same way as $\El(\bU)$,
 but allowing infinitely many elements of $\dX^*$ to occur in every representation. On the contrary, we
 always suppose that for every $i\in I$ the sum of all dimensions $n_x$ with $x\in\dE^*_i\cup\dF^*_i$ is finite.
 The last condition looks indispensable, since even when one considers
 the simplest case $\#I=1,\,\dE=\set x,\,\dF=\set y,\,x\sim y$ (which means square matrices under
 conjugation), the classification of representations of infinite dimension is a wild problem. 

 To deal with such infinite representations we first establish a general result concerning infinite matrices over
 bimodules.

 \begin{definition}\label{a21}
  Let $\bA$ be a locally finite dimensional category, $\bB$ be its full subcategory and $\bU$ be an $\bA$-bimodule. We
 say that $\bU$ is \emph{triangular with respect to $\bB$} if, for every indecomposable objects $A,B,C$, where $B,C\in\bB$
 and $A\notin \bB$, $\bA(C,A)\bU(B,C)=\bU(B,C)\bA(A,B)=0$.
 \end{definition} 
 
 The following lemma is obvious.
 \begin{lemma}\label{a22}
 Let $\bU$ be triangular with respect to $\bB$. For any object $A\in\bA$ choose a decomposition $A\iso A_1\+A_2$, where
 $A_1\in\bB$ and $A_2$ has no direct summands from $\bB$. For a morphism $a\in\bA(A,A')$ or an element $u\in\bU(A,A')$
 denote, respectively, by $a_1$ or $u_1$ its component from $\bA(A_1,A'_1)$ or $\bU(A_1,A'_1)$. If $a\in\bA(A,A')$ is a
 morphism in $\El(\bU)$ from $u\in\bU(A,A)$ to $v\in\bU(A',A')$ (i.e. $au=va$), then $a_1$ is a morphism from $u_1$ to $v_1$.
 Especially if $a$ is an isomorphism $u\to v$, then $a_1$ is an isomorphism $u_1\to v_1$.
 \end{lemma} 
 
 \begin{lemma}\label{a23}
 Let $\bA$ be a locally finite dimensional category that is a union of a chain $\bA_1\subseteq\bA_2\subseteq\bA_3\subseteq\dots$ of full
 subcategories. Suppose that $\bU$ is an $\bA$-bimodule that is triangular with respect to each $\bA_i$. Denote by $\bA^\8$
 the category of infinite direct sums $A=\bigoplus_{i=1}^\8 A_i$, where $A_i$ is an object from $\bA_i$ with no direct summands from
 $\bA_{i-1}$, and by $\bU^\8$ the natural extension of $\bU$ onto $\bA^\8$. For each element $u$ from $\bU^\8(A,A)$
 denote by $u_n$ its restriction onto $\bigoplus_{i=1}^na_i$. If $u,v$ are two elements such that $u_n\iso v_n$ for all $n$, then $u\iso v$. 
 \end{lemma} 
 {\sc Proof.}\ 
  First suppose the field $\fK$ \emph{uncountable}.
   Consider the sets of isomorphisms $\Iso(u_n,v_n)$ and the natural mappings $\pi_{mn}:\Iso(u_m,v_m)\to\Iso(u_n,v_n)\ (m>n)$
 arising from the triangularity condition. These sets can be considered as algebraic (even affine) varieties, then $\pi_{mn}$ are open
 morphisms of these varieties. In particular, the images $\im\pi_{m1}$ form a decreasing chain of non-empty open subsets in
 $\Iso(u_1,v_1)$. Hence their intersection is also non-empty (cf. for instance \cite{gab}). Take an element $a_1$ from this intersection
 and set $X_n=\pi_{n1}^{-1}(a_1)\ (n\ge2)$. Again they are algebraic varieties and $X'_n=\pi_{n2}(X_n)$ are non-empty
 open subsets of $X_2$, thus there is an element $a_2\in\bigcap_{n=2}^\8X'_n$. Continuing this process, we get a sequence $a_n$
 of elements from $\Iso(u_n,v_n)$ such that $\pi_{mn}(a_m)=a_n$ for all $m>n$. This sequence defines an isomorphism $a:u\to v$.

 If $\fK$ is arbitrary, take its uncountable extension $\ti\fK$ and consider extensions of $\bA$ and $\bU$ to
 $\ti\fK$. It is easy to see that $u\iso v$ \iff their extensions are isomorphic, which accomplishes the proof.\label{
} \hfill$\Box$ 

 Note that using Lemma \ref{12} one can obtain the following analogue of Proposition \ref{1a2} (with almost the same proof).

 \begin{proposition}\label{a230}
 We use the suppositions and notations of Lemma \ref{a23}. If $A=\bigoplus_{i=1}^\8A_i\in\bA^\8$ and $a\in\End A$, denote by
 $a_i$ the component of $a$ belongings to $\End A_i$. The category $\El(\bU^\8)$ is $\om$-local (in particular,
 Krull--Schmidt). Moreover, if
 $u\in\bU^\8(A,A)$ is an  indecomposable element from $\El(\bU^\8)$ and $a\in\End u$, then either all $a_i$ are invertible
 or all of them are nilpotent. 
 \end{proposition} 

\noindent
 Now we define \emph{infinite $\dX$-words} as sequences  $w=\dots x_1r_1x_2r_2x_3\dots r_{m-1}x_m \dots,$
 which are one-side or two-side infinite, subject to conditions (a) and (b) of Definition \ref{a12}\,(1) and such that for each
 $i$ the set $\setsuch{k}{x_k\in\dE_i\cup\dF_i}$ is finite. We apply to such words all terminology from Definitions
 \ref{a12}\,(2--4) and \ref{a13}\,(1,2,5ab). Then we can extend Theorem \ref{a14} to infinite representations.

  \begin{theorem}\label{a24}
  Isomorphism classes of indecomposable infinite representations of a bunch of chains are in \oc with equivalence
 classes of infinite string data. Moreover, every infinite representation uniquely decomposes into a direct sum of
 indecomposable ones. 
 \end{theorem} 
 {\sc Sketch of the proof}\ (more details will appear in \cite{Bond}).\ 
 Let $\dX_m$  be the bunch of chains with the index set $I_m=\{1,2,$ $\dots,m\}$, the same chains $\dE_i,\dF_i$
 and the same relation $\sim$, $\bA_m=\bA(\dX_m)$. Then we are in the situation of Lemma \ref{a23}.
 We define representations corresponding to infinite string data just as it has been done in \cite{mp,mp1} for finite case. 
 One can show  that all of them are indecomposable and their endomorphism rings are local. So we only have to prove 
 that there are no more indecomposable infinite representations.

  For each representation $M\in\El^\8(\bU)$ we denote by $M_m$ the
 restriction of $M$ onto $\dX_m$, given by all matrices $M_{xy}$ with $x,y\in\dX_m^*$.
 Lemma \ref{a23} implies that $M\iso N$ \iff $M_m\iso N_m$ for every $m$.
 Suppose that $M$ is infinite and indecomposable and
 consider an indecomposable direct summand $L$ of a representation $M_m$. The reduction procedure
 and the explicit description of strings and bands from \cite{mp1} immediately imply the following facts.
 \begin{enumerate}
 \item 
  $L$ cannot be a band or a bispecial string.
 \item
  If $L$ is a usual or a special string, there is an integer $m'>m$ and an indecomposable direct summand
 $L'$ of $M_{m'}$ such that the word $w$ from the string datum corresponding to $L$ is a part of the
 word $w'$ from the string datum corresponding to $L'$.
 \item
  If $K$ is another indecomposable direct summand of $M_m$, the number $m'>m$ and the representation
 $L'$ from (2) can be chosen common for $L$ and $K$.
 \end{enumerate}
 It implies the first statement of the theorem. The Krull--Schmidt property follows from Proposition \ref{a230}.
 \hfill$\Box$ 

\begin{center}
{\bf Acknowledgement}
\end{center}
{\small The authors would like to thank the referee for a very 
careful reading of the paper  and numerous remarks.  
The first author would like to thank V.M.~Bondarenko, \\
C.-M.~Rin\-gel and W.~S\"orgel for  helpful discussions of the results of this article.}

\end{document}